Lorenz Friess

Zusammenfassung

Sei G ein endlicher, ebener, zusammenhängender, schlichter Graph. Die Vielecke von G seien alle bis auf höchstens eins - etwa $V^*$ - Dreiecke, k bezeichne die Anzahl der Ecken von $V^*$. Dann gibt gibt es mindestens $2^{k-3}$ vierfarbige Färbungen der Ecken von G, derart, daß Ecken, die durch eine Kante verbunden sind, verschieden gefärbt sind. Diese Färbungen sind auf $V^*$ verschieden.

Abstract

Let G be a finite planar connected graph without loops or multiple edges. All minimal circuits except one - say $C^*$ - are triangles. Let k be the number of vertices of $C^*$. There are at least $2^{k-3}$ colorings of the vertices of G with four colors, vertices connected by edges colored differently. These colorings are different on $C^*$.

1. Definitionen, Bezeichnungen, Satz

Begriffe, die hier nicht erklärt werden, sind immer im Sinne von WAGNER Graphentheorie [1] zu verstehen. Dieses Buch wird mit [W, Seitenzahl] zitiert.

Im weiteren werden nur ebene, zusammenhängende, schlichte Graphen betrachtet.

In Analogie zum Grad $\gamma(E, G)$ einer Ecke $E \in G$ [W, 11] bezeichne $\gamma(V, G)$ die Anzahl der Kanten eines Weges $V \in G$ [W, 30]. Wenn innerhalb (oder außerhalb) eines geschlossenen Weges V (geschlossener Weg ist ein Kreis) keine Ecke und keine Kante von G liegt, so heißt V Vieleck.

$\gamma(G) = \max \gamma(V)$, V Vieleck von G, heißt der Grad von G. G heißt maximal, wenn $\gamma(G) = 3$, d.h. wenn alle Vielecke von G Dreiecke sind. ("maximal eben" in [W, 106])

Ein Graph G heißt fastmaximal, wenn für alle Vielecke V bis auf höchstens eins (etwa $V^*$) von G gilt $\gamma(V) = 3$. Die Ecken von $V^*$ heißen Randecken von G.

Für ein Rad R bezeichnet der Grad $\gamma(R)$ die Anzahl der Ecken auf dem Rand, also $\gamma(R) = |R| - 1$, die Ecke in der Mitte nennen wir Nabe. $R_i$ bezeichnet das Rad mit $\gamma(R) = i$.

Ein Graph G heißt innen-zusammenhängend (kurz i-Graph) wenn G - {Randecken} zusammenhängend und nicht leer ist (zusammenhängend [W, 32]). (Graph G - Ecke E ist der Graph G ohne E und ohne alle mit E inzidenten Kanten, [W, 22].) Räder sind i-Graphen.

Unter einer n-farbigen Färbung, $n \in N$ eines Graphen G verstehen wir eine Abbildung f der Menge der Ecken von G in { 1, 2, ... , n}. Eine Färbung heißt zulässig, wenn für die Ecken E, E' jeder Kante gilt: $f(E) \neq f(E')$.

4-farbige, zulässige Färbungen werden im folgenden kurz Färbungen genannt, andere Färbungen werden nicht betrachtet.

**Satz**
**Jeder fast-maximale Graph G besitzt mindestens $2^{\gamma(G)-3}$ Färbungen.**
**Diese unterscheiden sich auf dem Rand von G.**

## 2. Heuristische Überlegung

Versucht man einen maximalen Graphen G zu färben, so beginnt man etwa mit einem Dreieck und färbt dessen Ecken $E_1$, $E_2$, $E_3$ mit den Farben 1, 2 und 3. (Wir benutzen zur Bezeichnung der vier Farben die Ziffern 1, 2, 3 und 4, s. oben.) Mit dieser Festlegung entfällt der Faktor 24 für die Anzahl der Färbungen der von den Permutationen der vier Farben kommt.

Ausgehend von diesem Dreieck - $G_3$ - färben wir eine weitere Ecke, - $E_4$ - die mit mindestens zwei Ecken von $G_3$ verbunden ist. Wir färben weiter (zufällig oder systematisch) eine Ecke nach der anderen. Zweierlei fällt dabei auf:

1. bei vielen Ecken, vor allem am Anfang des Färbens, hat man mehrere Möglichkeiten

2. irgendwann kommt man i. a. in die Lage, daß man eine fünfte Farbe braucht. Letzteres kann man umgehen durch Umfärben einer oder mehrerer bereits gefärbter Ecken. Das bedeutet: man hat bei einer vorangehenden Ecke mit Wahlmöglichkeit "falsch" gewählt.

Wenn es gelingt sich alle Möglichkeiten (beim sukzessiven Färben) zu merken, so wird man - wenn es überhaupt eine (vierfarbige) Färbung des Graphen gibt - diese (und überhaupt alle) finden. Ein solches Vorgehen ist also optimal.

Wir betrachten als Beispiel den Graphen von Fig Heu1.

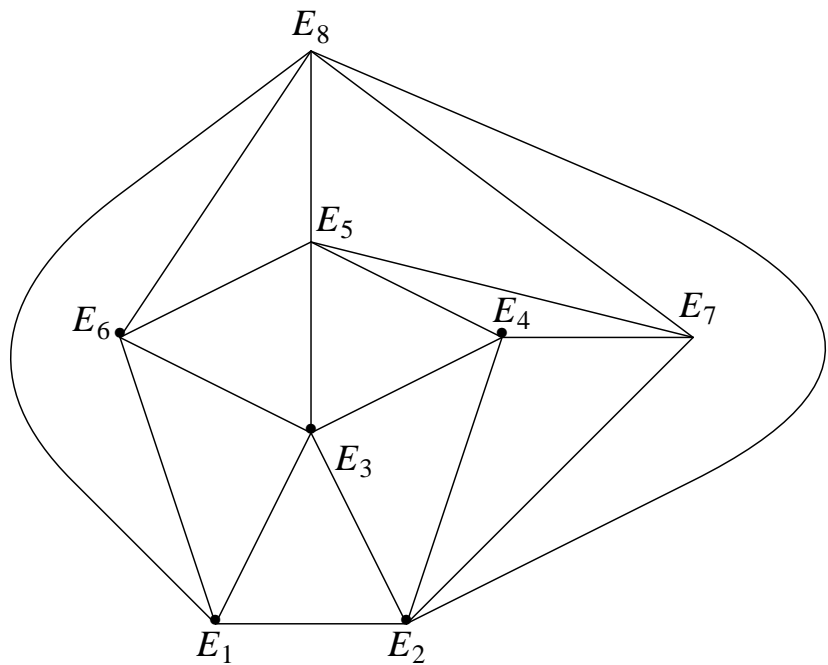

Fig Heu1

Wir beginnen die Färbung bei einem Dreieck, dem "Basisdreieck". Seine Ecken $E_1$, $E_2$, $E_3$ werden mit 1, 2, 3 gefärbt.

Nun kommt eine Ecke hinzu. Sie soll mit wenigstens zwei Kanten mit den Ecken des Basisdreiecks verbunden sein, Fig Heu2.

Es sind zwei Farben für $E_4$ möglich: 1 und 4. Wir denken uns den Rand des Graphen orientiert und notieren die Farben der Ecken des Randes in entsprechender Reihenfolge, also:

| $E_1$ | $E_2$ | $E_4$ | $E_3$ |
|-------|-------|-------|-------|
| 1     | 2     | 1     | 3     |
| 1     | 2     | 4     | 3     |

Die fünfte Ecke wählen wir so, daß wir den Teilgraphen Fig Heu3 haben.

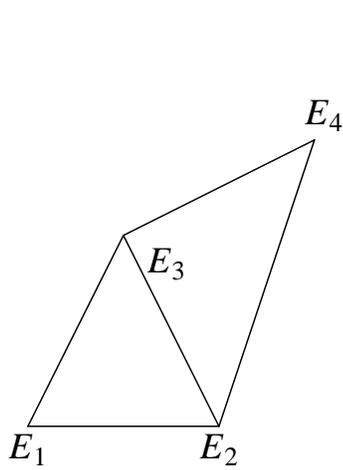
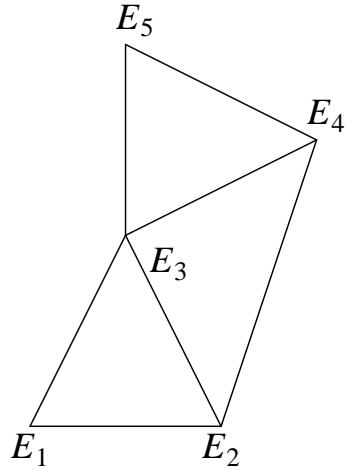
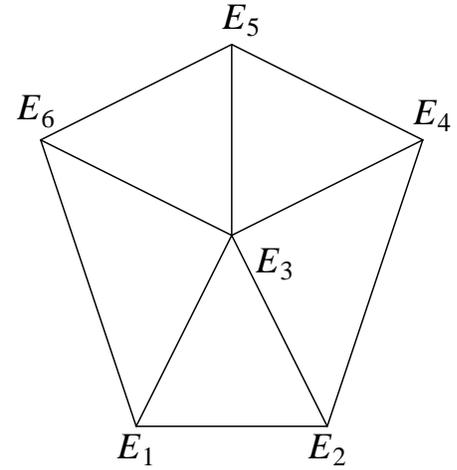

Fig Heu2                Fig Heu3                        Fig Heu4

Insgesamt sind vier Färbungen möglich:

| $E_1$ | $E_2$ | $E_4$ | $E_5$ | $E_3$ |
|-------|-------|-------|-------|-------|
| 1     | 2     | 1     | 2     | 3     |
| 1     | 2     | 1     | 4     | 3     |
| 1     | 2     | 4     | 1     | 3     |
| 1     | 2     | 4     | 2     | 3     |

Die Notierung der Färbungen ist mühsam und unübersichtlich, vereinfacht können wird schreiben für Fig Heu2 (für die möglichen Färbungen):

| $E_1$ | $E_2$ | $E_4$ | $E_3$ |
|-------|-------|-------|-------|
|       |       | 1     |       |
| 1     | 2     |       | 3     |
|       |       | 4     |       |

| $E_1$ | $E_2$ | $E_4$ | $E_5$ | $E_3$ |
| --- | --- | --- | --- | --- |
|  |  |  | 2 |  |
|  |  | 1 |  |  |
|  |  |  | 4 |  |
| 1 | 2 |  |  | 3 |
|  |  |  | 2 |  |
|  |  | 4 |  |  |
|  |  |  | 1 |  |

Es dürfte klar sein wie eine solche Tabelle - wir nennen sie "Färbungsschema" (FS) - zu lesen ist. Es wird sich im weiteren zeigen, daß diese Art die möglichen Färbungen zu beschreiben mehr als nur eine Arbeitserleichterung ist.

Wir setzen die Färbung des Beispiels fort: für den Graph Fig Heu4 ergibt sich das Färbungsschema (daneben alle Färbungen als "normale" Liste)

| $E_1$ | $E_2$ | $E_4$ | $E_5$ | $E_6$ |  | $E_1$ | $E_2$ | $E_4$ | $E_5$ | $E_6$ |
| --- | --- | --- | --- | --- | --- | --- | --- | --- | --- | --- |
|  |  |  | 2 | 4 |  | 1 | 2 | 1 | 2 | 4 |
|  |  | 1 |  |  |  |  |  |  |  |  |
|  |  |  | 4 | 2 |  | 1 | 2 | 1 | 4 | 2 |
| 1 | 2 |  |  |  | = |  |  |  |  |  |
|  |  |  | 2 | 4 |  | 1 | 2 | 4 | 2 | 4 |
|  |  | 4 |  |  |  |  |  |  |  |  |
|  |  |  | 1 | 2, 4 |  | 1 | 2 | 4 | 1 | 2 |
|  |  |  |  |  |  | 1 | 2 | 4 | 1 | 4 |

Die Färbungen der Ecke $E_3$, die jetzt im Inneren liegt, treten nicht im Färbungsschema auf. Das Färbungsschema beschreibt Färbungen durch Angabe der Farben der Ecken des Randes. Zu diesen Färbungen existieren "Fortsetzungen" nach innen. Wenn wir die (zu verschiedenen Färbungen des Randes möglicherweise verschiedenen) Fortsetzungen nach innen nicht aufschreiben, "vergessen" wir Information, die für die weitere Konstruktion der Färbungen offenbar überflüssig ist. (Mit diesen inneren Ecken werden keine weiteren Ecken mehr durch Kanten verbunden und daher bedeuten ihre Farben keine Bedingung für die weiteren Ecken.)

Bemerkung: G wird als plättbar vorausgesetzt. Hier wird deutlich, daß diese Voraussetzung wesentlich ist.

| $E_1$ | $E_2$ | $E_7$ | $E_5$ | $E_6$ | $E_1$ | $E_2$ | $E_7$ | $E_5$ | $E_6$ |
|---|---|---|---|---|---|---|---|---|---|
| | | 4, 3 | 2 | 4 | 1 | 2 | 4 | 2 | 4 |
| | | | | | 1 | 2 | 3 | 2 | 4 |
| | | 3 | 4 | 2 | 1 | 2 | 4 | 2 | 4 |
| 1 | 2 | | | | | | | | |
| | | 1, 3 | 2 | 4 | 1 | 2 | 1 | 2 | 4 |
| | | | | | 1 | 2 | 3 | 2 | 4 |
| | | 3 | 1 | 2, 4 | 1 | 2 | 3 | 1 | 2 |
| | | | | | 1 | 2 | 3 | 1 | 4 |

Das FS beschreibt 4 + 3 Färbungen! Davon sind jedoch zwei gleich.

Die Ecke $E_4$ fällt nach innen, es gibt keine Spalte von $E_4$ im FS.

Das FS wird komplizierter. Wir lösen das daraus resultierende Problem folgendermaßen: wir schließen nicht eine Ecke (hier $E_7$) und danach noch eine Ecke (hier $E_8$) an, sondern das Rad um $E_7$ mit den Ecken $E_2, E_4, E_5, E_8$, mit $f(E_7) = 3$, dessen FS (des Rades) wir ja bereits kennen, dann die Kante $E_8 E_6$, dann die Kante $E_8 E_1$.

Diese Vorgehensweise wird unten als Anschluß Rad und Anschluß Kante beschrieben.

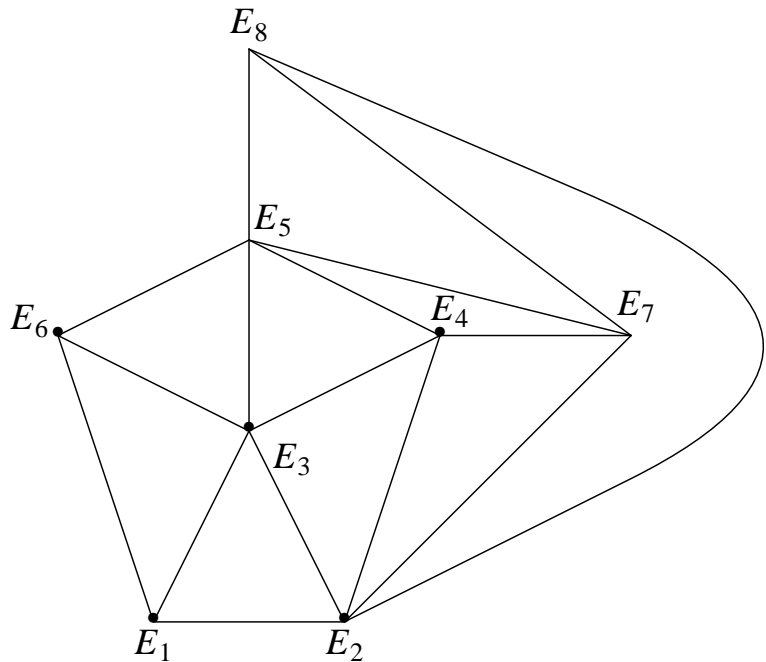

Fig Heu5

Als FS haben wir (Fig Heu5):

| $E_1$ | $E_2$ | $E_8$ | $E_5$ | $E_6$ |
|-------|-------|-------|-------|-------|
|       |       |       | 2     | 4     |
|       |       | 1     |       |       |
|       |       |       | 4     | 2     |
|       |       |       |       |       |
| 1     | 2     |       |       |       |
|       |       |       |       |       |
|       |       |       | 2     | 4     |
|       |       | 4     |       |       |
|       |       |       | 1     | 2, 4  |

Mit der Kante $E_8 E_6$ ($E_5$ fällt nach innen):

| $E_1$ | $E_2$ | $E_8$ | $E_6$ |
|-------|-------|-------|-------|
|       |       | 1     | 4, 2  |
| 1     | 2     |       |       |
|       |       | 4     | 2     |

Mit der Kante $E_8 E_1$ ($E_6$ fällt nach innen):

| $E_1$ | $E_2$ | $E_8$ |
|-------|-------|-------|
| 1     | 2     | 4     |

## 3. Weitere Bezeichnungen

### Länge einer Ecke E

Zur Vereinfachung der Ausdrucksweise wird mit der Länge einer Ecke E auf dem Rand eines fm Graph die Länge der Spalte von E im FS bezeichnet.

### Anschluß einer Ecke E

Beim Anschluß einer Ecke E an einen Graph G sei ihr relativer Grad >= 2, d.h. sie wird mit (mindestens) zwei Kanten mit G verbunden. (Relativ deshalb, weil der Grad von E durch weiter Anschlüsse noch größer werden kann.)

### Anschluß Rad, Anschluß fm Graph

Unter dem Anschluß eines Rades R an einen Graphen G verstehen wir die Identifikation von Randecken des Rades mit Randecken von G wie es aus den Beispielen genau hervorgeht. Wir schreiben dafür z. B. H = G + R.

Der Anschluß eines fm Graphen G' an einen fm Graphen G bedeutet ganz analog die Identifikationen von Ecken (nicht unbedingt alle Ecken!) der Ränder der beiden Graphen, Schreibweise H = G + G'. (Wieviele Ecken der Ränder von G, G' geht später aus dem Zusaamhang hervor.)

Wichtig ist, daß die Randecken von beiden Graphen "lückenlos" identifiziert werden, dann ist nämlich der Graph, der durch diese Anschlüsse entsteht, fast-maximal.

Zu einem fm Graph G denken wir uns ein Dreieck als Basisdreieck gewählt, das ist ein Dreieck, von dem genau eine Seite auf dem Rand von G liegt. Das FS beginnt mit diesem Dreieck wie in Abschnitt 2 erläutert. D. h. die Ecken der Basisseite haben die Länge 1, die weiteren Ecken sind dann jeweils doppelt so lang bis auf die zwei letzten (auf dem Rand), die gleich lang sind.

4. Definitionen gewisser Klassen von Graphen

$K_1 = \{$ G | G ist Rad, $| G | \geq 4$ $\}$

$K_2 = \{$ G | $G_1 +$ Rad, $G_1 \in K_1$, G ist ein fm i-Graph $\}$

$K_n = \{$ G | $G_{n-1} +$ Rad, $G_{n-1} \in K_{n-1}$, G ist ein fm i-Graph $\}$

Bemerkung: Damit $G = G_n +$ Rad ein fm i-Graph ist, wird das Rad mit mindestens 3 Ecken an $G_n$ angeschlossen. (siehe auch Abschnitt 11.)

$K_n^* = \{$ $G = G_n +$ "erlaubter Anschluß", $G_n \in K_n$ $\}$

Erlaubter Anschluß ist Anschluß eines Bogens gemäß Abschnitt 8, Anschluß einer Ecke gemäß Abschnitt 9, Anschluß einer Kante gemäß Abschnitt 10.
(Auch mehrere Anschlüsse an verschiedenen Stellen des Randes von G, nicht überschneidend.)

Sei $K = \bigcup_n K_n^*$

Wir beweisen zunächst den Satz für $G \in K$ und betrachten danach fm Graphen $G \notin K$.

Färbungsschema für $R_5$, FS($R_5$), siehe Fig Heu4

| $E_1$ | $E_2$ | $E_4$ | $E_5$ | $E_6$ | Färbungen |
|---|---|---|---|---|---|
| | | | 2 | 4 | 1 2 1 2 4 |
| | | 1 | | | |
| | | | 4 | 2 | 1 2 1 4 2 |
| 1 | 2 | | | | |
| | | | 2 | 4 | 1 2 4 2 4 |
| | | 4 | | | |
| | | | 1 | 2, 4 | 1 2 4 1 2 und 1 2 4 1 4 |

Für $E_3$ gibt es keine Spalte im FS da diese Ecke innen liegt, d.h. nicht auf dem Rand von G.

Es folgen die FS für $R_6$ und für $R_7$.

FS($R_6$)

| $E_1$ | $E_2$ | $E_4$ | $E_5$ | $E_6$ | $E_7$ | | | | | |
|---|---|---|---|---|---|---|---|---|---|---|
| | | | | 1 | 4, 2 | | | | | |
| | | | 2 | | | | | | | |
| | | | | 4 | 2 | | | | | |
| | | 1 | | | | | | | | |
| | | | | 1 | 2, 4 | | | | | |
| | | | 4 | | | | | | | |
| | | | | 2 | 4 | | | | | |
| 1 | 2 | | | | | | | | | |
| | | | | 4 | 2 | | 1 | 2 | | 1 | 4 |
| | | | 2 | | | 2 | | | 4 | |
| | | | | 1 | 4, 2 | | 4 | 1, 2 | | 2 | 1, 4 |
| | | 4 | | | | = 1 | | | = 1 | |
| | | | | 4 | 2 | | 1 | 2 | | 1 | 4 |
| | | | 1 | | | 4 | | | 2 | |
| | | | | 2 | 4 | | 2 | 1 | | 4 | 1 |

Hinter "=" steht der Block nach Permutation(1, 4), dann (2, 4) für spätere Erklärung.

Die Spalte "Block" gehört nicht zum FS sondern wird zur Erläuterung verwendet.
(Es ist $f(E_1) = 1$, $f(E_2) = 2$.)

**Block I**

| $E_4$ | $E_5$ | $E_6$ | $E_7$ | $E_8$ | Block | $E_6$ | $E_7$ | $E_8$ | $E_6$ | $E_7$ | $E_8$ |
|---|---|---|---|---|---|---|---|---|---|---|---|
|  |  |  | 2 | 4 |  |  | 1 | 4 |  | 2 | 1 |
|  |  | 1 |  |  |  | 2 |  |  | 4 |  |  |
|  |  |  | 4 | 2 |  |  | 4 | 1 |  | 1 | 2 |
|  | 2 |  |  |  | I = 1 |  | = 2 |  |  |  |  |
|  |  |  | 2 | 4 |  |  | 1 | 4 |  | 2 | 1 |
|  |  | 4 |  |  |  | 4 |  |  | 1 |  |  |
|  |  |  | 1 | 2,4 |  |  | 2 | 1,4 |  | 4 | 2,1 |

| $E_4$ |
|---|
| 1 |

**Block II**

| $E_4$ | $E_5$ | $E_6$ | $E_7$ | $E_8$ | Block |
|---|---|---|---|---|---|
|  |  |  | 4 | 2 |  |
|  |  | 1 |  |  |  |
|  |  |  | 2 | 4 |  |
|  | 4 |  |  |  | II = I |
|  |  |  | 4 | 2 |  |
|  |  | 2 |  |  |  |
|  |  |  | 1 | 4, 2 |  |

**Block III**

| $E_4$ | $E_5$ | $E_6$ | $E_7$ | $E_8$ | Block |
|---|---|---|---|---|---|
|  |  |  | 2 | 4 |  |
|  |  | 4 |  |  |  |
|  |  |  | 1 | 2, 4 |  |
|  | 2 |  |  |  | III = I |
|  |  |  | 2 | 4 |  |
|  |  | 1 |  |  |  |
|  |  |  | 4 | 2 |  |

| $E_4$ |
|---|
| 4 |

**Block IV**

| $E_4$ | $E_5$ | $E_6$ | $E_7$ | $E_8$ | Block | $E_6$ | $E_7$ | $E_8$ | $E_6$ | $E_7$ | $E_8$ |
|---|---|---|---|---|---|---|---|---|---|---|---|
|  |  |  | 1 | 2,4 |  |  | 2 | 1,4 |  | 4 | 1,2 |
|  |  | 4 |  |  |  | 4 |  |  | 2 |  |  |
|  |  |  | 2 | 4 |  |  | 1 | 4 |  | 1 | 2 |
|  | 1 |  |  |  | IV = 2 |  | = 4 |  |  |  |  |
|  |  |  | 1 | 4,2 |  |  | 2 | 4,1 |  | 4 | 2,1 |
|  |  | 2 |  |  |  | 1 |  |  | 1 |  |  |
|  |  |  | 4 | 2 |  |  | 4 | 1 |  | 2 | 1 |

Es gibt also (nach Permutation, Zeilen vertauschen) nur zwei Blocktypen, I und IV, für gekoppelte Spalten.

Wenn wir in den längsten Spalten, also zu $E_7$ und $E_8$, die jeweils paarweise auf naheliegende Weise zusammengehörenden Farben betrachten ("2-Block"), stellen wir fest, daß es nur folgende "Typen" gibt:

|   | 2 | 4 |   |   | 1 | 4, 2 |   |   | 2 | 4 |
|---|---|---|---|---|---|------|---|---|---|---|
| 1 |   |   | , | 2 |   |      | , | 4 |   |   |
|   | 4 | 2 |   |   | 4 | 2    |   |   | 1 | 2, 4 |

Nach Permutation der Farben (und geeigneter Reihenfolge der Zeilen) also nur zwei Typen.

Wir stellen fest: in den letzten drei Spalten eines FS eines Rades treten nur Blöcke auf der Form:

|   | 2 | 4 |   |   | 1 | 4, 2 |
|---|---|---|---|---|---|------|
| 1 |   |   | , | 2 |   |      |
|   | 4 | 2 |   |   | 4 | 2    |

Von den (analog gebildeten) vier 4-Blöcken, Block I bis IV im FS oben, sind die ersten drei gleich (nach Permutation, geeignete Reihenfolge der Zeilen), d.h. es gibt zwei Typen von 4-Blöcken. Diese Blöcke finden wir auch schon im $\mathrm{FS}(R_6)$.

Als Beispiel zeigen wir Block III = Block II:

|   | 2 | 4 |   |   |   | 4 | 2    |   |   |   | 4 | 2 |
|---|---|---|---|---|---|---|------|---|---|---|---|---|
| 4 |   |   |   |   | 2 |   |      |   | 1 |   |   |   |
|   | 1 | 2, 4 |   |   |   | 1 | 4, 2 |   |   |   | 2 | 4 |

2           = 4           = 4

|   | 2 | 4 |   |   |   | 4 | 2 |   |   |   | 4 | 2    |
|---|---|---|---|---|---|---|---|---|---|---|---|------|
| 1 |   |   |   |   | 1 |   |   |   |   | 2 |   |      |
|   | 4 | 2 |   |   |   | 2 | 4 |   |   |   | 1 | 4, 2 |

Für die Räder $R_8$, $R_9$, ... gibt es offenbar jeweils ganz analog ein FS. Es kommen auch nur dieselben Blöcke vor. (Das ist leicht einzusehen: die obere Hälfte von $\mathrm{FS}(R_6)$ ist gleich Block IV, die untere Hälfte von $\mathrm{FS}(R_6)$ ist gleich Block I, durch eine weitere Ecke ($E_8$, $\mathrm{FS}(R_7)$) erhalten wir zu Block IV zweimal Block I, zu Block I erhalten wir Block I und Block IV. Genauso liefert eine weitere Ecke ($E_9$, $\mathrm{FS}(R_8)$) wieder entsprechende Blöcke im $\mathrm{FS}(R_8)$.)

Die FS können wir auch mit der anderen Orientierung bilden, d.h. die Ecken des Randes in umgekehrter Reihenfolge färben. Außerdem können wir das Basisdreieck beliebig wählen. (Wegen der Symmetrie des Rades trivial.)

Die Anzahl der verschiedenen Färbungen (auf dem Rand verschieden, es existiert eine Fortsetzung nach innen) beträgt $2^{\gamma(G)-3}$ (mindestens).

Bemerkung: Die Anzahl der Färbungen ist $\geq 2^{\gamma(G)-3}$. Diese zusätzlichen Färbungen haben wir da, wo bei einem Element eine zweite Farbe eingetragen ist. Diese Färbungen werden bei der Anzahlangabe nicht mehr erwähnt.

Anfügen einer Kante von $E_i$ nach $E_{i+2}$, i= 4, ..., $\gamma(G) - 2$ liefert ein FS mit halber Länge und im übrigen wieder mit denselben Eigenschaften. Dazu betrachten wir als Beispiel FS($R_6$) und verbinden $E_4$ mit $E_6$. Färbungen, für die die Ecken $E_4$ und $E_6$ die gleiche Farbe haben, sind unzulässig und entfallen, als FS haben wir:

| $E_1$ | $E_2$ | $E_4$ | $E_5$ | $E_6$ |
|---|---|---|---|---|
|  |  |  | 2 | 4 |
|  |  | 1 |  |  |
|  |  |  | 4 | 2 |
| 1 | 2 |  |  |  |
|  |  |  | 1 | 4, 2 |
|  |  | 4 |  |  |
|  |  |  | 2 | 4 |

also gleich FS($R_5$).

Der Anschluß einer solchen Kante ist trivial (im Sinne des Färbungsproblems), da eine innere Ecke vom Grad 3 entsteht. Dies dient nur als Illustration der Eigenschaften von FS. Wenn wir - als weiteres Beispiel - eine Kante von $E_5$ nach $E_7$ anfügen erhalten wir als FS:

| $E_1$ | $E_2$ | $E_4$ | $E_5$ | $E_7$ | (Perm) | $E_1$ | $E_2$ | $E_4$ | $E_5$ | $E_7$ |
|---|---|---|---|---|---|---|---|---|---|---|
|  |  |  | 2 | 4 |  |  |  |  | 4 | 2 |
|  |  | 1 |  |  |  |  |  | 1 |  |  |
|  |  |  | 4 | 2 |  |  |  |  | 2 | 4 |
| 1 | 2 |  |  |  | (2, 4) | 1 | 4 |  |  |  |
|  |  |  | 2 | 4 |  |  |  |  | 4 | 2 |
|  |  | 4 |  |  |  |  |  | 2 |  |  |
|  |  |  | 1 | 2, 4 |  |  |  |  | 1 | 4, 2 |

also FS($R_5$). Dazu noch die Kante von $E_2$ nach $E_5$ gibt das FS

| $E_1$ | $E_2$ | $E_5$ | $E_7$ |
|---|---|---|---|
|  |  | 4 | 2 |
| 1 | 2 |  |  |
|  |  | 1 | 4, 2 |

Das ist gleich dem FS von $R_4$.

Da man sich ein FS aus Blöcken - wie oben dargestellt - zusammengesetzt denken kann, braucht man sich für das Anfügen einer Kante von einer Ecke zur übernächsten Ecke nur anzusehen was mit den Färbungen, die durch diese Blöcke beschrieben werden, passiert: die Färbungen können wieder durch die Blöcke beschrieben werden.

Es dürfte klar sein, was unter zusammengehörenden Blöcken zu verstehen ist: zwei Blöcke gehören zusammen, wenn sie in der vorherigen Spalte zu demselben Element gehören.

Wir betrachten nochmal $FS(R_7)$ und den Fall, daß eine Kante von $E_5$ nach $E_7$ angefügt wird. Als FS haben wir dann:

| $E_1$ | $E_2$ | $E_4$ | $E_5$ | $E_7$ | $E_8$ | (Perm) | $E_4$ | $E_5$ | $E_7$ | $E_8$ |
|---|---|---|---|---|---|---|---|---|---|---|
| | | | | 4 | 2 | | | | | |
| | | | 2 | | | | | | | |
| | | | | 1 | 2, 4 | | | | | |
| | | 1 | | | | | | | | |
| | | | | 2 | 4 | | | | | |
| | | | 4 | | | | | | | |
| | | | | 1 | 4, 2 | | | | | |
| 1 | 2 | | | | | | | | | |
| | | | | 1 | 4, 2 | | | | 2 | 4, 1 |
| | | | 2 | | | | | 1 | | |
| | | | | 4 | 2 | | | | 4 | 1 |
| | | 4 | | | | (1, 2) | 4 | | | |
| | | | | 2 | 4 | | | | 1 | 4 |
| | | | 1 | | | | | 2 | | |
| | | | | 4 | 2 | | | | 4 | 1 |

Aus Typ I und Typ I wird Typ IV, aus Typ I und Typ IV wird Typ I.

Typ IV tritt nur zusammen mit Typ I auf, d. h. nur in FS mit entsprechend großer Anzahl an Färbungen. (s. Abschnitt 6)

Zu jeder Ecke des Randes eines fm Graphen gehört eine Spalte im FS. Die Reihenfolge der Spalten entspricht der Reihenfolge der Ecken auf dem Rand des Graphen. Die Ecken der Basisseite haben die Länge 1. Es gibt eine dritte Spalte der Länge 1 ($\gamma(G) = 3$) oder die beiden letzten Spalten sind gleich lang und haben die Länge $2^{\gamma(G)-3}$ (mit 2 Farben auf einigen Positionen, vgl. Blöcke für gekoppelte Ecken). Diese beiden Ecken bzw. Spalten nennen wir gekoppelte Ecken bzw. gekoppelte Spalten.

$FS_2$

In einem FS gibt es nur 4-Blöcke der folgenden Typen (vgl. Abschnitt 5) (Für $\gamma(G) = 4$ nur 2-Block.)

|  |  |  |  |  |  |  |  |  |
|---|---|---|---|---|---|---|---|---|
|  | 1 |  |  | 1 | 4 |  | 1 | 2, 4 |
| 2 |  |  | 2 |  |  | 4 |  |  |
|  | 4 |  |  | 4 | 1 |  | 2 | 4 |
| 1 |  |  | 1 |  |  | 1 |  |  |
|  | 1 |  |  | 1 | 4 |  | 1 | 4, 2 |
| 4 |  |  | 4 |  |  | 2 |  |  |
|  | 2 |  |  | 2 | 1, 4 |  | 4 | 2 |
|  | Bl1 |  |  | Bl2 |  |  | Bl3 |  |

$FS_3$

Zu jedem Element in einer Spalte gehört auf naheliegende Weise ein Teil-FS ("Element" ist selbsterklärend). Dies wird an nachstehendem FS erläutert.

| $E_1$ | $E_2$ | $E_3$ | $E_4$ |
|---|---|---|---|
|  |  |  | 2 |
|  |  | 1 |  |
|  |  |  | 4 |
|  | 2 |  |  |
|  |  |  | 2 |
|  |  | 4 |  |
|  |  |  | 1 |
| 1 |  |  |  |
|  |  |  | 4 |
|  |  | 1 |  |
|  |  |  | 2 |
|  | 4 |  |  |
|  |  |  | 4 |
|  |  | 2 |  |
|  |  |  | 1 |

Zu "2" in der Spalte von $E_2$ gehört

| $E_2$ | $E_3$ | $E_4$ |
|---|---|---|
|  |  | 2 |
|  | 1 |  |
|  |  | 4 |
| 2 |  |  |
|  |  | 2 |
|  | 4 |  |
|  |  | 1 |

und zu "4" in Spalte $E_3$ gehört

| $E_3$ | $E_4$ |
|---|---|
|  | 2 |
| 4 |  |
|  | 1 |

jeweils mit den entsprechenden Teilen aller weiterer Spalten - soweit vorhanden.

Das FS oben besteht aus Bl1 für die Spalten von $E_1, E_2, E_3$, dann zweimal Bl2 für $E_2, E_3, E_4$.

Zum Ende dieses Abschnitts eine Übersicht, alle Formen von Typ I und Typ IV die man durch Permutation der Farben erhält, jeweils nur eine der möglichen Reihenfolgen der Zeilen.

Typ I

| Nr. | | | | | Nr | | | | |
|---|---|---|---|---|---|---|---|---|---|
|  |  |  | 1 | 4 |  |  |  | 1 | 2 |
|  |  | 2 |  |  |  |  | 4 |  |  |
|  |  |  | 4 | 1 |  |  |  | 2 | 1 |
| 1 | 1 |  |  |  | 2 | 1 |  |  |  |
|  |  |  | 1 | 4 |  |  |  | 1 | 2 |
|  |  | 4 |  |  |  |  | 2 |  |  |
|  |  |  | 2 | 1, 4 |  |  |  | 4 | 1, 2 |
|  |  |  | 2 | 4 |  |  |  | 2 | 1 |
|  |  | 1 |  |  |  |  | 4 |  |  |
|  |  |  | 4 | 2 |  |  |  | 1 | 2 |
| 3 | 2 |  |  |  | 4 | 2 |  |  |  |
|  |  |  | 2 | 4 |  |  |  | 2 | 1 |
|  |  | 4 |  |  |  |  | 1 |  |  |
|  |  |  | 1 | 2, 4 |  |  |  | 4 | 2, 1 |

- 15 -

```
        1                        2
           2    4                   1    4
5    4                   6    4
           4    2                   4    1
        2                    1
           1    4,2                  2    4,1
```

Typ IV

Nr.                         Nr

```
           1    4,2                  1    2,4
        2                         4
           4    2                     2    4
1    1                   2    1
           1    2,4                  1    4,2
        4                         2
           2    4                     4    2

           2    4,1                  2    1,4
        1                         4
           4    1                     1    4
3    2                   4    2
           2    1,4                  2    4,1
        4                         1
           1    4                     4    1

           4    2,1                  4    1,2
        1                         2
           2    1                     1    2
5    4                   6    4
           4    1,2                  4    2,1
        2                         1
           1    2                     2    1
```

Nach Vertauschen von Zeilen sind 1=2, 3=4, 5=6.

Die Blöcke Typ I, 1 bis 6 werden mit B1, ..., B6 bezeichnet.
Die Blöcke Typ IV, 1 bis 6 werden mit C1, ..., C6 bezeichnet. (Siehe Abschnitt 5.)

Blöcke):

```
        A1              A2              A3
            2               1               1
        1               2               4
            4               4               2
    2               1               1
            2               1               1
        4               4               2
            1               2               4
1               2               4
            4               4               2
        1               2               4
            2               1               1
    4               4               2
            4               4               2
        2               1               1
            1               2               4
```

Zur Erläuterung  die Blöcke in umgekehrter Orientierung:

```
        A1              A2              A3
2               1               1
    1               2               4
4               4               2
        2               1               1
2               1               1
    4               4               2
1               2               4
            1               2               4
4               4               2
    1               2               4
2               1               1
        4               4               2
4               4               2
    2               1               1
1               2               4
```

Bemerkung: aus der Struktur der Blöcke ergibt sich, wie ein FS aufgebaut sein kann.

Zur Schreibweise: aus typographischen Gründen schreiben wir manchmal A1 statt $A_1$, A2 statt $A_2$ usw. auch für B1, ..., C1, ...

Die Blöcke $A_i'$ sind beschreibbar durch

$$
A_1' \quad\quad A_2' \quad\quad A_3'
$$

$$
\begin{array}{ccc}
1 & 2 & 4 \\
2 & 1 & 1 \\
4 & 4 & 2 \\
\end{array}
$$

$$
1 \quad , 2 \quad , 4
$$

$$
\begin{array}{ccc}
1 & 2 & 4 \\
4 & 4 & 2 \\
2 & 1 & 1 \\
\end{array}
$$

dies soll so dargestellt werden:

$$
A_1 = \begin{matrix} A_2' \\ A_1' \\ A_3' \end{matrix} \quad\quad
A_2 = \begin{matrix} A_1' \\ A_2' \\ A_3' \end{matrix} \quad\quad
A_3 = \begin{matrix} A_1' \\ A_3' \\ A_2' \end{matrix}
$$

Die Blöcke $B_i$, $C_i$ kommen nur am Ende eines FS vor, haben also keine Nachfolger. Zu $A_1$ sind als Nachfolger möglich:

$$
A_1 \begin{matrix} B_3, B_4, C_3 \\ \\ B_5, B_6, C_5 \end{matrix}
$$

Tatsächlich gibt es aber nicht alle neun Kombinationen. Die Kopplung der Spalten in $B_i$, $C_i$ entsteht ja dadurch, daß eine Farbe "verboten" ist (durch eine Kante zu einer Ecke der Länge 1), diese verbotene Farbe ist für beide Blöcke dieselbe, also gibt es nur:

$$
A_1 \begin{matrix} B_3 \\ B_5 \end{matrix} \quad , A_1 \begin{matrix} B_4 \\ C_5 \end{matrix} \quad , A_1 \begin{matrix} C_3 \\ B_6 \end{matrix}
$$

Für $A_2$, $A_3$ findet man entsprechend:

$$
A_2 \begin{matrix} B_1 \\ C_1 \end{matrix} \quad , A_2 \begin{matrix} B_2 \\ B_4 \end{matrix} \quad , A_2 \begin{matrix} C_1 \\ B_3 \end{matrix}
$$

$$
A_3 \begin{matrix} B_1 \\ C_3 \end{matrix} \quad , A_3 \begin{matrix} B_2 \\ B_4 \end{matrix} \quad , A_3 \begin{matrix} C_1 \\ B_3 \end{matrix}
$$

Für $A_1$, $A_2$, $A_3$

$$
A_1 \begin{matrix} A_2 \\ A_3 \end{matrix} \quad , A_2 \begin{matrix} A_3 \\ A_1 \end{matrix} \quad , A_3 \begin{matrix} A_1 \\ A_2 \end{matrix}
$$

Graphen R (Ecken vom Rand) identifiziert, so ist dies zulässig nur möglich, wenn es passende - d. h. gleiche - Färbungen gibt. Ob dies der Fall ist hängt ab von den Blöcken, die die Färbungen dieser Ecken beschreiben. Die Ränder von G und R sind entgegengesetzt orientiert.

Liste der Kombinationen von Blöcken mit passender Färbung - soweit vorhanden - als Beispiel. Bei jeder Kombination hat man sich die beiden Blöcke mit entgegengesetzter Orientierung zu denken:

Tabelle 4-4

|    | A1     | A2     | A3     |
|----|--------|--------|--------|
| A1 | 1 2 4 1 | 2 1 4 1 | 4 2 1 2 |
| A2 |        | 2 1 4 2 | 2 4 2 4 |
| A3 |        |        | 4 2 1 4 |

|    | B1     | B2     | B3     | B4     | B5     | B6     |
|----|--------|--------|--------|--------|--------|--------|
| A1 | 1 4 2 1 | 1 2 4 1 | -      | 2 1 2 1 | -      | 4 1 4 1 |
| A2 | -      | 1 4 1 2 | 2 4 1 2 | 2 4 1 2 | 4 1 4 2 | -      |
| A3 | 1 4 1 4 | -      | 2 1 2 4 | -      | 4 2 1 4 | 4 2 1 4 |
| B1 | 1 4 2 1 |        |        |        |        |        |
| B2 | 1 4 2 1 | 1 2 4 1 |        |        |        |        |
| B3 | -      | -      | 2 1 4 2 |        |        |        |
| B4 | -      | 2 4 2 1 | 2 1 4 2 | 2 4 1 2 |        |        |
| B5 | -      | -      | 4 1 4 2 | -      | 4 2 4 2 |        |
| B6 | 4 2 4 1 | -      | -      | -      | 4 1 2 4 | 4 1 2 4 |

|    | C1     | C3     | C5     |
|----|--------|--------|--------|
| A1 | -      | 2 1 2 1 | 4 1 4 1 |
| A2 | 1 2 1 2 | -      | 4 2 4 2 |
| A3 | 1 4 1 4 | 2 4 2 4 | -      |
| B1 | -      | -      | 1 2 1 4 |
| B2 | -      | 1 4 1 2 | -      |
| B3 | -      | -      | 2 1 2 4 |
| B4 | 2 4 2 1 | -      | -      |
| B5 | -      | 4 2 1 2- |        |
| B6 | 1 4 2 4 | -      | -      |
| C1 | -      |        |        |
| C3 | 1 2 1 2 | -      |        |
| C5 | 1 4 1 4 | 2 1 2 4 | -      |

Kombinationen für die "-" eingetragen ist, können nicht zur Konstruktion des FS(G+R) verwendet werden, da keine passenden Färbungen für die Ecken existieren.

Tabelle der Blöcke und passende Farbungen für drei gemeinsame Spalten.

<div align="center">Tabelle 3-3</div>

|     | A1    | A2    | A3    |       |       |       |
|-----|-------|-------|-------|-------|-------|-------|
| A1  | 2 1 2 |       |       |       |       |       |
| A2  | 4 2 4 | 1 2 1 |       |       |       |       |
| A3  | 2 4 2 | 1 2 1 | 1 4 1 |       |       |       |

|     | B1    | B2    | B3    | B4    | B5    | B6    |
|-----|-------|-------|-------|-------|-------|-------|
| A1  | 2 1 4 | 2 1 2 | 4 1 2 | 4 1 2 | 1 4 2 | 1 2 4 |
| A2  | 4 2 1 | 4 2 1 | 1 2 4 | 4 2 1 | 1 2 4 | 2 4 1 |
| A3  | 2 4 1 | 4 2 1 | 1 4 2 | 4 2 1 | 1 4 2 | 2 4 1 |
| B1  | 4 1 4 | 2 1 4 | 1 4 2 | 2 1 4 | 2 4 1 | 4 2 1 |
| B2  |       | 2 1 2 | 1 2 4 | 1 4 2 | 1 4 2 | 4 2 1 |
| B3  |       |       | 4 2 4 | 1 2 4 | 4 1 2 | 1 4 2 |
| B4  |       |       |       | 1 2 1 | 4 2 1 | 4 2 1 |
| B5  |       |       |       |       | 2 4 2 | 1 4 2 |
| B6  |       |       |       |       |       | 1 4 1 |

|     | B1    | B2    | B3    | B4    | B5    | B6    |
|-----|-------|-------|-------|-------|-------|-------|
| C1  | 4 1 4 | 2 1 2 | 4 2 4 | 4 1 2 | 2 1 2 | 2 1 4 |
| C3  | 4 2 1 | 4 2 1 | 1 2 4 | 1 2 1 | 1 2 4 | 1 2 4 |
| C5  | 2 4 1 | 2 1 2 | 1 4 2 | 1 2 1 | 1 2 4 | 1 2 4 |

|     | C1    | C3    | C5    |
|-----|-------|-------|-------|
| A1  | 2 1 2 | 1 2 4 | 2 1 2 |
| A2  | 4 2 4 | 1 2 1 | 1 2 1 |
| A3  | 2 1 2 | 1 4 1 | 2 1 2 |

|     | C1    | C3    | C5    |
|-----|-------|-------|-------|
| C1  | 2 1 2 | 4 2 4 | 2 1 2 |
| C3  |       | 1 2 1 | 1 2 1 |
| C5  |       |       | 1 2 1 |

Genauso erhält man noch eine Tabelle für Blöcke mit zwei gemeinsamen Spalten.

In den folgenden Abschnitten ist G - soweit nicht anders definiert - ein fm Graph mit FS.

In diesem Abschnitt betrachten wir den Anschluß eines Rades R an einen fm Graph G mit FS. Mit Anschluß R(n, k) bezeichnen wir den Anschluß eines Rades R mit $\gamma(R) = n$, wobei k Ecken von R mit Randecken von G identifiziert werden. Es sind nur R(n, k) mit k ≥ 3 von Interesse (sonst ist G + R kein i-Graph).

Für die Fälle R(n=2k-2, k) gibt es eine spezielle besonders einfache Konstruktion des FS(G+R), s. unten Fall m = k. In diesen Fällen ist $\gamma(G) = \gamma(G + R)$.

Gegeben sei ein fm Graph G mit Färbungsschema FS. Den Graph G stelle man sich in FigRad "unter" $E_{g_1}$ , $E_{g_2}$ , ... $E_{g_{k-1}}$ , $E_{g_k}$ vor, $E_1, E_n$ sind durch eine Kante verbunden.

An G wird ein Rad R angeschlossen. Behauptung: es gibt ein FS für G + R.

Beweis durch Konstruktion des FS(G+R).

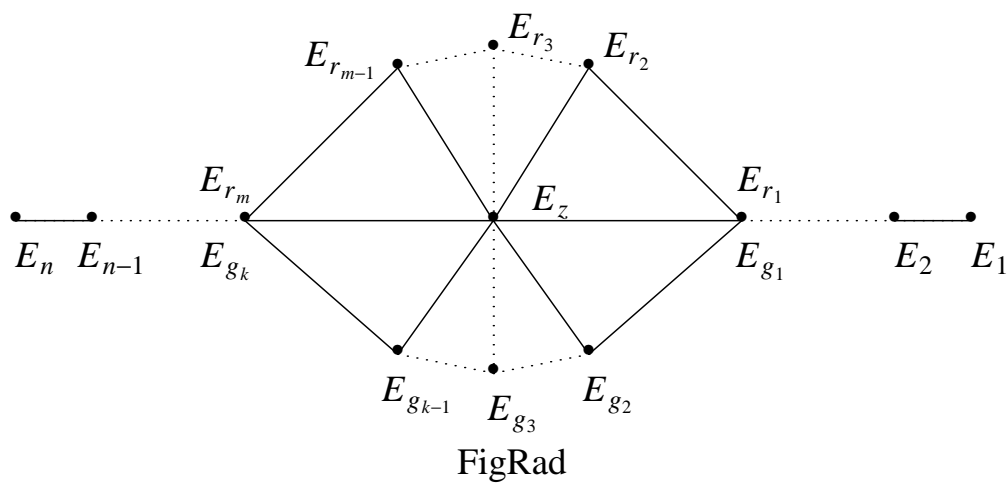

FigRad

Bezeichnungen:

$E_1$ , ..., $E_n$ die Ecken vom Rand von G, $E_1 E_2$ Basisseite von G, $E_n E_{n-1}$ gekoppelte Spalten von G; $E_{r_1}$ , ..., $E_{r_{m+k}}$ die Ecken vom Rand von R, $E_{r_2}, \ldots, E_{r_{m-1}}$ die Ecken von R, die nicht mit Ecken von G identifiziert werden, $E_{r_1}$ , ..., $E_{r_m}$ die Ecken von R, die auf dem Rand von G + R liegen, also $\gamma(R) = m + k - 2$. $E_{g_1}$ , ..., $E_{g_k}$ bezeichnen die Ecken von G, die mit Ecken von R identifiziert werden, nämlich $E_{r_1}$ mit $E_{g_1}$ , $E_{r_{m+i}}$ mit $E_{g_{k-i}}$ , i=2,...,k-2.

R wird mit $k \geq 3$ Ecken an G angeschlossen (für $k \leq 2$ ist der entstehende Graph kein i-Graph). $m \geq 2$, m = 2 bedeutet Kante von $E_{r_1}$ nach $E_{r_2}$. Es gilt: $\gamma(G + R) = \gamma(G) + (m-2) - (k-2) = \gamma(G) + m - k$.

Bemerkung zum Fall m = k.

Wir betrachten zuerst kurz den speziellen Fall m = k. Er ist besonders einfach - dies ist für den weiteren Beweis jedoch ohne Bedeutung.

In diesem Fall wählen wir $f(E_z) = 3$ und $f(E_{r_i}) = f(E_{g_i})$, $i = 2, \ldots, k-1$, d.h. wir haben wieder dasselbe FS: FS(G+R)=FS(G).

Bevor wir andere Wertepaare für m, k betrachten geben wir eine Übersicht der Beweisschritte:

Gegeben ist der fm Graph G mit FS(G), anzuschließen ist das Rad R.

Das FS für den fm Graph G+R erhalten wir so:

Die Basisseite von G+R sei die von G und für die Randecken $E_1$ bis $E_{g_1}$ wählen wir das FS von G, d. h. den entsprechenden Teil des FS von G.

Wir wählen $E_{g_1}$ und $E_{g_2}$ als Ecken der Basisseite von R (die Basisseite ist frei wählbar).

Für die Ecken von R, die auf dem Rand von G + R liegen, (die "neuen" Randecken außer $E_{g_1}$, $E_{g_k}$) wählen wir das TeilFS bestehend aus den Spalten zu $E_{r_1}, \ldots, E_{r_m}$ vom FS(R) und setzen damit das FS(G+R) zu **jedem** Element von $E_{g_1}$ fort.

Für die Ecken $E_{g_k}$ und darüber verwenden wir wiederum Teilschemata vom FS(G) (evtl. mehrfach benutzt: trotzdem sind die Färbungen im FS von G+R alle verschieden: sie unterscheiden sich auf $E_{r_2}, \cdots E_{r_{m-1}}$.)

Mit jeder Ecke von R, die auf dem Rand von G + R liegt, wird die Anzahl der Färbungen verdoppelt, danach mit jeder Ecke von G die nach innen fällt die Anzahl der Färbungen halbiert. (Die genaue Anzahl liefert das FS.)

Zu zeigen: das Rad kann "verträglich" an G angeschlossen werden, d. h. es gibt in FS(G) und in FS(R) für die Ecken, die identifiziert werden, gleiche ("passende") Farben. Dies geschieht unter Verwendung der Blöcke aus denen man sich die FS von G und R aufgebaut denken kann.

Als Spezialfälle muß auf die Lage der gekoppelten Spalten eingegangen werden, z. B. wenn die gekoppelten Spalten von G nach innen fallen, d. h. mit Ecken von R identifiziert werden, oder wenn k=3, dann liegen die gekoppelten Spalten von R beide auf dem Rand von G+R.

Die Konstruktion des FS von G + R im einzelnen.

Bezeichnung: $E_x < E_y$ soll bedeuten $E_x$ liegt im Sinne der Orientierung des Randes von G (oder R) vor $E_y$.

Der Anfang der Konstruktion ist immer gleich: es wird für $E_1$ bis $E_{g_1}$ das FS von G als (Anfang des) FS von G+R übernommen.

Wir unterscheiden bei der weiteren Konstruktion die Lage der gekoppelten Spalten, da für diese die Blöcke anders aufgebaut sind als im übrigen FS.

Es folgt eine Reihe von Fallunterscheidungen.

Zuerst als Beispiel der Fall k > 4, m > 4.

Die gekoppelten Ecken $E_{r_{m+k-3}}$, $E_{r_{m+k-4}}$ von R werden identifiziert mit $E_{g_3}$, $E_{g_4}$.

Sei zunächst $E_{n-1} > E_{g_k}$.

Die Basisseite von G + R sei die von G. Für die Basisseite und die Randecken bis $E_{g_1}$ wählen wir das FS von G. Für jedes Element der Spalte von $E_{g_1}$ und den beiden dazugehörenden Elementen der Spalte von $E_{g_2}$ haben wir ein FS von R mit $E_{g_1}$, $E_{g_2}$ ($= E_{r_1}$, $E_{r_{k+m-2}}$) als Basisseite.

Für das folgende betrachten wir $f(E_{g_1}) = 1$, dann ist $f(E_{g_2}) = 2$ oder 4, für $f(E_{g_1}) = 2$ und $f(E_{g_1}) = 4$ verlaufen alle Überlegungen analog. Mit diesem FS(R) (Teil-FS bis zur Spalte von $E_{r_{m-1}}$) wird das FS von G + R durch Anhängen an jedes Element von $E_{g_1}$ definiert. Die so definierten Färbungen für G + R (zunächst für die Randecken bis einschließlich $E_{r_{m-1}}$) bilden ein FS für den bis jetzt betrachteten Teil des Randes von G+R.

Haben nun die Färbungen von R Fortsetzungen für $E_{r_m}$, ..., $E_{r_{m+k-3}}$($= E_{g_3}$) derart, daß sie zu den Färbungen von G auf den entsprechenden Ecken passen? (Die möglichen Färbungen auf den Ecken von G liegen fest weil wir das FS von G verwenden für den Anfang des FS von G + R.)

Wir sehen uns die Spalte von $E_{r_{m-1}}$ (im FS von R) an (einen Block, die Argumentation ist für alle Blöcke gleich) diese Spalte ist die letzte, die einfach in das FS(G+R) übernommen ("abgeschrieben") wurde.

Für $E_{r_{m-3}}$, $E_{r_{m-2}}$, $E_{r_{m-1}}$, $E_{r_m}$ liegt Block $A_1$, $A_2$ oder $A_3$ vor. Es kann nur einer dieser drei Fälle vorliegen, welcher es tatsächlich ist, hängt von m ab.

Wir betrachten andererseits einen Block aus dem FS von G, und zwar zu einem Element der Spalte von $E_{g_{k-2}}$.

Er kann - in Abhängigkeit von $f(E_{g_{k-2}})$ - ebenfalls $A_1$, $A_2$ oder $A_3$ sein, und zwar gehört zu $f(E_{g_{k-2}}) = 1$ bzw. $f(E_{g_{k-2}}) = 2$ bzw. $f(E_{g_{k-2}}) = 4$ Block $A_1$ bzw. $A_2$ bzw. $A_3$. Welcher tatsächlich vorliegt hängt von k ab.

Zu den Farben in der Spalte von $E_{r_{m-1}}$ können wir jetzt die Färbung des Randes (die Konstruktion des FS von G + R) fortsetzen indem wir Färbungen des Randes von G verwenden:

Zu $f(E_{r_{m-1}}) = 1$ sind einerseits vom (FS vom) Rad $f(E_{r_m}) = 2$ und $f(E_{r_m}) = 4$ möglich und andererseits - für alle Blöcke $A_1$, $A_2$, $A_3$ aus dem FS von G - sind $f(E_{g_k}) = 2$ und $f(E_{g_k}) = 4$ möglich, und zwar mit den zugehörigen Fortsetzungen vom FS von G (für die Ecken oberhalb $E_{g_k}$) bis $E_n$ und diese Färbungen unterscheiden sich auf dem Rand (von G + R). Dieselben Fortsetzungen werden mehrfach verwendet!

Zu $f(E_{r_{m-1}}) = 2$ sind einerseits vom FS von R $f(E_{r_m}) = 1$ und $f(E_{r_m}) = 4$ möglich und andererseits - für alle Blöcke $A_1$, $A_2$, $A_3$ aus dem FS von G - sind $f(E_{g_k}) = 1$ und $f(E_{g_k}) = 4$ möglich. (Es kann mehrere Möglichkeiten geben, z. b. $f(E_{g_k}) = 1$ zweimal. Es ist egal welche "1" mit ihren Fortsetzungen verwendet wird, die Fortsetzungen sind gleich, vgl. Abschnitt 5.)

Zu $f(E_{r_{m-1}}) = 4$ sind $f(E_{r_m}) = 1$ und $f(E_{r_m}) = 2$ vom FS von R möglich und andererseits $f(E_{g_k}) = 1$ und $f(E_{g_k}) = 2$ möglich für alle Blöcke vom FS von G.

Damit R so angeschlossen werden kann und das FS für G + R konstruiert werden kann, müssen wir noch zeigen, daß es im FS von R und im FS von G für die Ecken, die jeweils identifiziert werden, passende (d. h. gleiche) Färbungen gibt.

$f(E_{r_m}) = f(E_{g_k})$ haben wir gerade betrachtet, jetzt ist noch $f(E_{r_{m+1}}) = f(E_{g_{k-1}})$, ..., $f(E_{r_{k+m-3}}) = f(E_{g_3})$, $f(E_{r_{k+m-2}}) = f(E_{g_2})$, zu zeigen. D. h. zu dem Vorgänger von $E_{g_k}$, das ist $E_{g_{k-1}}$ - dessen Farbe steht fest - muß es unter den Fortsetzungen von $E_{r_m}$ auf $E_{r_{m+1}}$ eine passende (d.h. gleiche) Farbe geben usw., also zu zeigen: es gibt Färbungen sodaß $f(E_{r_{m+1}}) = f(E_{g_{k-1}})$, ..., $f(E_{r_{k+m-3}}) = f(E_{g_3})$, $f(E_{r_{k+m-2}}) = f(E_{g_2})$,

Für diese Ecken wird - wie bereits festgestellt - das FS von R durch die Blöcke $A_1$, $A_2$, $A_3$ beschreiben.

Dieselben Blöcke haben wir haben wir im FS von G. Es müssen also zusammenpassen die Blöcke $A_1$, $A_2$, $A_3$ mit einem der drei $A_1$, $A_2$, $A_3$ in umgekehrter Orientierung. Es gibt immer zueinander passende Färbungen gleichgültig welcher Block von G und welcher Block von R zusammenpassen müssen wie man aus der Abschnitt 6, Tabelle 4-4 ersehen kann.

Wir stellen den Sachverhalt nochmal anders dar: haben wir z. B. $f(E_{g_{k-i}}) = f(E_{r_{m+i}})$, i=1, ..., k-4, so ist in FS(G) $f(E_{g_{k-i-1}}) = 2$ **oder** 4 und in FS(G) $f(E_{r_{m+i+1}}) = 2$ **und** 4. Daher gibt es zwei passende Färbungen (zwei gleiche, je eine in FS(G) und in FS(R)).

Diese Konstruktion ist so möglich Spalte für Spalte bis die gekoppelten Spalten von R anzuschließen sind, d.h. für $E_{r_i}$, $i < m + k - 4$.

Wir betrachten nun die möglichen Blöcke des FS von R für die gekoppelten Spalten. Die gekoppelten Spalten gehören zu den Ecken $E_{r_{m+k-4}}$, $E_{r_{m+k-3}}$.

Welcher Block für R vorliegt hängt ab von $f(E_{g_2})$. ($E_{g_2}$ war als eine Ecke der Basisseite von R gewählt und damit geht $f(E_{g_2})$ in das FS(R) ein.) Es kann $f(E_{g_2}) = 2$ oder $f(E_{g_2}) = 4$ sein. (Wir betrachten $f(E_{g_1}) = 1$.) Damit erhalten wir zwei FS für R. (Zu einem Element von $E_{g_1} = E_{r_{k+m-1}}$ haben wir zwei Elemente in (der Spalte von) $E_{g_2} = E_{r_{k+m-2}}$ und für diese zwei verschiedene FS für R.) Wenn $f(E_{g_2}) = 2$ haben wir im FS von R einen der Blöcke B1, B6, C3. Die entsprechenden Ecken werden mit $E_{g_3}$, ... identifiert. Für diese Ecken der Block in FS(G) ist A2. Es gibt immer passende Färbungen, Abschnitt 6, Tabelle 3-3.

|  | $E_{r_{m+k-5}}$ | $E_{r_{m+k-4}}$ | $E_{r_{m+k-3}}$ | $E_{r_{m+k-2}}$ |
|---|---|---|---|---|
|  | $E_{g_5}$ | $E_{g_4}$ | $E_{g_3}$ | $E_{g_2}$ |
|  |  | 2 | 4, 1 |  |
|  | 1 |  |  |  |
|  |  | 4 | 1 |  |
| R, C3 |  |  |  | 2 |
|  |  | 2 | 1, 4 |  |
|  | 4 |  |  |  |
|  |  | 1 | 4 |  |
|  | 1 |  |  |  |
|  |  | 2 |  |  |
|  | 4 |  |  |  |
|  |  |  | 1 |  |
|  | 1 |  |  |  |
|  |  | 4 |  |  |
|  | 2 |  |  |  |
| G, A2 |  |  |  | 2 |
|  | 4 |  |  |  |
|  |  | 2 |  |  |
|  | 1 |  |  |  |
|  |  |  | 4 |  |
|  | 4 |  |  |  |
|  |  | 1 |  |  |
|  | 2 |  |  |  |

Wir stellen fest, daß es zu jeder Färbung von $E_{r_{m+k-5}}$, $E_{r_{m+k-4}}$, $E_{r_{m+k-3}}$, $E_{r_{m+k-2}}$ passende Färbungen der entsprechenden Ecken von G gibt.

| $E_{g_5}$ | $E_{g_4}$ | $E_{g_3}$ | $E_{g_2}$ |
|---|---|---|---|
| 1 | 2 | 1 | 2 |
| 1 | 2 | 4 | 2 |
| 1 | 4 | 1 | 2 |
| 4 | 2 | 1 | 2 |
| 4 | 2 | 4 | 2 |
| 4 | 1 | 4 | 2 |

| | $E_{r_{m+k-5}}$ | $E_{r_{m+k-4}}$ | $E_{r_{m+k-3}}$ | $E_{r_{m+k-2}}$ |
| --- | --- | --- | --- | --- |
| | $E_{g_5}$ | $E_{g_4}$ | $E_{g_3}$ | $E_{g_2}$ |
| | | 1 | 4 | |
| | 2 | | | |
| | | 4 | 1 | |
| R, B1 | | | | 2 |
| | | 1 | 4 | |
| | 4 | | | |
| | | 2 | 1, 4 | |

| | $E_{g_5}$ | $E_{g_4}$ | $E_{g_3}$ | $E_{g_2}$ |
| --- | --- | --- | --- | --- |
| | 1 | | | |
| | | 2 | | |
| | 4 | | | |
| | | | 1 | |
| | 1 | | | |
| | | 4 | | |
| | 2 | | | |
| G, A2 | | | | 2 |
| | 4 | | | |
| | | 2 | | |
| | 1 | | | |
| | | | 4 | |
| | 4 | | | |
| | | 1 | | |
| | 2 | | | |

die passenden Färbungen:

| | | | |
| --- | --- | --- | --- |
| 2 | 1 | 4 | 2 |
| 2 | 4 | 1 | 2 |
| 4 | 1 | 4 | 2 |
| 4 | 2 | 1 | 2 |
| 4 | 2 | 4 | 2 |

Wir stellen fest, daß es zu jeder Färbung von $E_{r_{m+k-5}}$, $E_{r_{m+k-4}}$, $E_{r_{m+k-3}}$, $E_{r_{m+k-2}}$ passende Färbungen der entsprechenden Ecken von G gibt.

| | $B_{r_{m+k-5}}$ $E_{g_5}$ | $B_{r_{m+k-4}}$ $E_{g_4}$ | $B_{r_{m+k-3}}$ $E_{g_3}$ | $B_{r_{m+k-2}}$ $E_{g_2}$ |
|---|---|---|---|---|
| | | 4 | 1 | |
| | 1 | | | |
| | | 2 | 4, 1 | |
| R, B6 | | | | 2 |
| | | 4 | 1 | |
| | 2 | | | |
| | | 1 | 4 | |
| | 1 | | | |
| | | 2 | | |
| | 4 | | | |
| | | | 1 | |
| | 1 | | | |
| | | 4 | | |
| | 2 | | | |
| G, A2 | 4 | | | 2 |
| | | 2 | | |
| | 1 | | | |
| | | | 4 | |
| | 4 | | | |
| | | 1 | | |
| | 2 | | | |

die passenden Färbungen:

| 1 | 4 | 1 | 2 |
|---|---|---|---|
| 1 | 2 | 4 | 2 |
| 1 | 2 | 1 | 2 |
| 2 | 4 | 1 | 2 |
| 2 | 1 | 4 | 2 |

Für $f(E_{g_2}) = 4$ gibt entsprechendes. Blöcke für R sind B2, B4, C3, für G Block A3.

| $E_{r_{m+k-5}}$ $E_{g_5}$ | $E_{r_{m+k-4}}$ $E_{g_4}$ | $E_{r_{m+k-3}}$ $E_{g_3}$ | $E_{r_{m+k-2}}$ $E_{g_2}$ |
|---|---|---|---|
|  | 1 | 2 |  |
| 4 |  |  |  |
|  | 2 | 1 |  |
| R, B2 |  |  | 4 |
|  | 1 | 2 |  |
| 2 |  |  |  |
|  | 4 | 1, 2 |  |
| 1 |  |  |  |
|  | 4 |  |  |
| 2 |  |  |  |
|  |  | 1 |  |
| 1 |  |  |  |
|  | 2 |  |  |
| 4 |  |  |  |
| G, A3 |  |  | 4 |
| 2 |  |  |  |
|  | 4 |  |  |
| 1 |  |  |  |
|  |  | 2 |  |
| 2 |  |  |  |
|  | 1 |  |  |
| 4 |  |  |  |

die passenden Färbungen:

| | | | |
|---|---|---|---|
| 4 | 1 | 2 | 4 |
| 4 | 2 | 1 | 4 |
| 2 | 1 | 2 | 4 |
| 2 | 4 | 1 | 4 |
| 2 | 4 | 2 | 4 |

| $E_{r_{m+k-5}}$ | $E_{r_{m+k-4}}$ | $E_{r_{m+k-3}}$ | $E_{r_{m+k-2}}$ |
| $E_{g_5}$ | $E_{g_4}$ | $E_{g_3}$ | $E_{g_2}$ |
|---|---|---|---|
|  | 1 | 2 |  |
| 4 |  |  |  |
|  | 2 | 1 |  |
| R, B4 |  |  | 4 |
|  | 2 | 1 |  |
| 1 |  |  |  |
|  | 4 | 2, 1 |  |
| 1 |  |  |  |
|  | 4 |  |  |
| 2 |  |  |  |
|  |  | 1 |  |
| 1 |  |  |  |
|  | 2 |  |  |
| 4 |  |  |  |
| G, A3 |  |  | 4 |
| 2 |  |  |  |
|  | 4 |  |  |
| 1 |  |  |  |
|  |  | 2 |  |
| 2 |  |  |  |
|  | 1 |  |  |
| 4 |  |  |  |

die passenden Färbungen:

| | | | |
|---|---|---|---|
| 4 | 1 | 2 | 4 |
| 4 | 2 | 1 | 4 |
| 1 | 2 | 1 | 4 |
| 1 | 4 | 2 | 4 |
| 1 | 4 | 1 | 4 |

| | $B_{r_{m+k-5}}$ | $B_{r_{m+k-4}}$ | $B_{r_{m+k-3}}$ | $B_{r_{m+k-2}}$ |
|---|---|---|---|---|
| | $E_{g_5}$ | $E_{g_4}$ | $E_{g_3}$ | $E_{g_2}$ |
| | | 4 | 2, 1 | |
| | 1 | | | |
| | | 2 | 1 | |
| R, C5 | | | | 4 |
| | | 4 | 1, 2 | |
| | 2 | | | |
| | | 1 | 2 | |
| | 1 | | | |
| | | 4 | | |
| | 2 | | | |
| | | | 1 | |
| | 1 | | | |
| | | 2 | | |
| | 4 | | | |
| G, A3 | | | | 4 |
| | 2 | | | |
| | | 4 | | |
| | 1 | | | |
| | | | 2 | |
| | 2 | | | |
| | | 1 | | |
| | 4 | | | |

die passenden Färbungen:

| | | | |
|---|---|---|---|
| 1 | 4 | 2 | 4 |
| 1 | 4 | 1 | 4 |
| 1 | 2 | 1 | 4 |
| 2 | 4 | 1 | 4 |
| 2 | 4 | 2 | 4 |
| 2 | 1 | 2 | 4 |

Für $f(E_{g_1}) = 4$ haben wir im FS(R) einen der Blöcke B2, B4, C5.

|  |  | Blöcke FS(R) | Block FS(G) |
|---|---|---|---|
| $f(E_{g_1}) = 1$ | $f(E_{g_2}) = 2$ | B1, B6, C3 | A2 |
|  | $f(E_{g_2}) = 4$ | B2, B4, C5 | A3 |
| $f(E_{g_1}) = 2$ | $f(E_{g_2}) = 1$ | B3, B5, C1 | A1 |
|  | $f(E_{g_2}) = 4$ | B2, B4, C5 | A3 |
| $f(E_{g_1}) = 4$ | $f(E_{g_2}) = 1$ | B3, B5, C1 | A1 |
|  | $f(E_{g_1}) = 2$ | B1, B6, C3 | A2 |

Wie bereits oben ("Übersicht der Beweisschritte") bemerkt, wird das FS(G+R) für $E_{g_k}$ bis $E_n$ durch die entsprechenden TeilFS(G) zu jedem Element von $E_{g_k}$ fortgesetzt. FS(G+R) hat dann offensichtlich alle Eigenschaften eines FS.

Bemerkung: Bei der Konstruktion des FS von G + R werden nicht immer alle möglichen Färbungen berücksichtigt. Dies entspricht dem Sachverhalt, daß es maximale Graphen gibt mit mehr als einer Färbung. (Hier hat ein maximaler Graph ja nur eine Färbung.)
Wenn man diese nicht berücksichtigten Färbungen untersucht, erhält man Bedingungen dafür, daß ein maximaler Graph mehr als eine Färbung besitzt.

Die Bedeutung der FS und deren Aufbau durch Blöcke ist nun deutlich.

Es folgt die Konstruktion des FS(G+R) mit zahlreichen Fallunterscheidungen. Diese Konstruktionen verwenden als Eigenschaft von R nur, daß es ein fm Graph mit FS ist (mit beliebig wählbarer Basisseite und Orientierung). Der Graph R muß kein Rad sein.

Wir unterscheiden die Konstruktion des FS(G+R) nach der Lage der gekoppelten Ecken von G und R zueinander und ob sie auf dem Rand von G+R liegen.

Es sei erinnert an $E_{g_1} = E_{r_1}$, $\gamma(R) = m + k - 2$, $E_{r_i} = E_{g_{m+k-i}}$, i=m+1, ..., m+k-2.

Mit $E_{g_{k+1}}$ (nicht $E_{g_{k+1}}$) wir die nächste Ecke auf dem Rand von G bezeichnet. In diesem Abschnitt sei $m \geq 3$. Die Fälle mit m=2 werden extra untersucht.

In den Abbildungen ist die Basisseite verstärkt dargestellt, die Kante zwischen gekoppelten Ecken ist durch eine zusätzliche Strichelung dargestellt. Für $E_{g_1} = E_2$ ist $f(E_{g_1}) = 2$, dann ist $f(E_{g_2}) = 1$ oder 4, die Konstruktion des FS(G+R) kann genauso durchgeführt werden. Die Fälle $E_{g_1} > E_2$, $E_{g_1} = E_2$, müssen erst bei $f(E_{g_n}) = f(E_{g_k})$ getrennt betrachtet werden. Kanten , die gestrichelt sind, können fehlen, dann sind die beiden Ecken natürlich identisch.

1. $E_{n-1} \geq E_{g_k}$, d.h. die gekoppelten Ecken von G liegen auf dem Rand von G+R.

1.1.

$k \geq 4$, $m \geq 4$, d. h. die gekoppelten Ecken von R liegen nicht auf dem Rand von G+R.

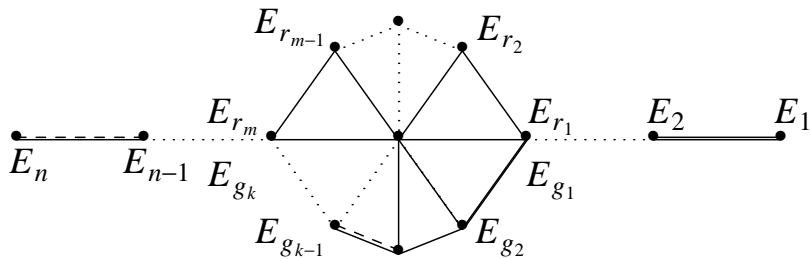

$E_{g_1} > E_2$, dieser Fall ist oben ausführlich dargestellt.

Das FS(G) zu den Ecken $E_{g_1}, ... E_{g_k}$ wird beschrieben durch normale Blöcke, also durch A1, A2 oder A3.

Das FS(R) zu den Ecken $E_{r_1}$ bis $E_{r_{m+k-4}}$ d. h. bis zu den gekoppelten Ecken $E_{r_{m+k-4}}$, $E_{r_{m+k-3}}$ wird beschrieben durch normale Blöcke. Also müssen passen $A_i$, $A_j$, i, j = 1,2,3. Dies ist immer erfüllt.

Der Block für die Ecken $E_{m+k-4}$, $E_{m+k-3}$, die gekoppelten Ecken von R, hängt ab von $f(E_{g_2})$. Sei $f(E_{g_1}) = 1$, dann ist $f(E_{g_2}) = 2$ oder 4 und der Block für G für die Ecken $E_{g_1}$ $E_{g_2}$ $E_{g_3}$ $E_{g_4}$ ist vom Typ A1.

Für $f(E_{g_2}) = 2$ haben wir für R den Block B1 oder B6 oder C3, da $f(E_{g_3} = E_{r_{m+k-3}}) \neq 2$ sein muß.

Den Tabellen entnehmen wir, daß es passenden Färbungen gibt. (Alle Kombinationen von A1, A2, A3 mit B1, B6, C3.)

Für $f(E_{g_2}) = 4$ haben wir für R den Block B2 oder B4 oder C5.

Den Tabellen entnehmen wir, daß es passenden Färbungen gibt.

Für $E_{g_1} = E_2$ ist dieselbe Konstruktion möglich, es ist dabei nur $f(E_{g_1}) = 2$ ($= f(E_2)$) zu beachten.

1.2.1 m=3. Wir setzen $f(E_{r_2}) = f(E_{g_2})$.

1.2.2 m=4, $E_{g_1} \geq E_2$

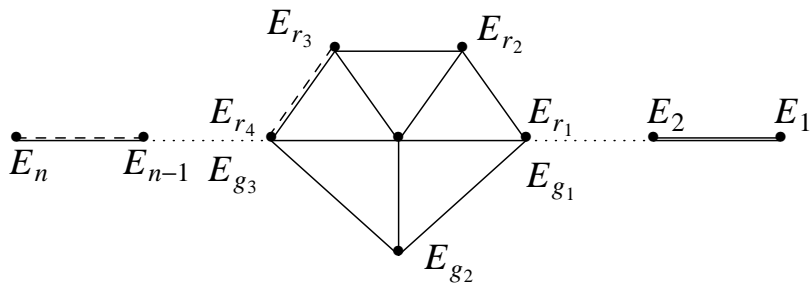

Wir betrachten die Lage von $E_{g_1}$, d. h. ob $E_{g_1} = E_2$ oder $E_{g_1} \neq E_2$:

$E_{g_1} = E_2 \Rightarrow f(E_{g_1}) = 2 \Rightarrow f(E_{g_2}) = 1$ oder $4$

$E_{g_1} > E_2$

$E_{g_1} = E_3 \Rightarrow f(E_{g_1}) = 1$ bzw. $4$ und $f(E_{g_2}) = 2, 4$ bzw. $f(E_{g_2}) = 1, 2$.

$E_{g_1} > E_3 \Rightarrow f(E_{g_1}) = 1, 2$ oder $1, 4$ oder $2, 4$.

In jedem dieser Fälle gibt es für einen Wert von $f(E_{g_1})$ zwei Werte für $f(E_{g_2})$ und daher mit $E_{g_1}, E_{g_2}$ als Basisseite zwei FS für R.

Damit wir ein FS für G+R erhalten, muß bei der Konstruktion der Block des FS(R), der die Kopplung $E_{r_3}, E_{r_4}$ beschreibt, durch einen Block ohne Kopplung ersetzt werden (in einem FS sind nur die beiden letzten Ecken gekoppelt, in G+R also $E_n, E_{n-1}$ und für diese Ecken haben wir die Kopplung des FS(G+R) per Konstruktion).

Im FS(R) sind für $E_{r_1}, ..., E_{r_4}$ abhängig von $f(E_{g_1})$, $f(E_{g_2})$ nur bestimmte Blöcke möglich:

| $f(E_{g_1})$ | $f(E_{g_2})$ | Block im FS(G) | Blöcke im FS(R) | |
|---|---|---|---|---|
| | | | m = 4 | m > 4 |
| 2 | 1 | A2 | B3 | B3, B5, C1 |
| | 4 | | B4 | B4, B2, C5 |
| | | | | |
| 1 | 2 | A1 | B1 | B1, B6, C3 |
| | 4 | | B2 | B2, B4, C5 |
| | | | | |
| 4 | 1 | A3 | B5 | B5, B3, C1 |
| | 2 | | B6 | B6, B1, C3 |

Wir betrachten die möglichen Färbungen von $E_{g_1}$.

$f(E_{g_1}) = 1$. Für $f(E_{g_2}) = 2$ haben wir im FS(R) Block B2, für $f(E_{g_2}) = 4$ haben wir im FS(R) Block B3.

Aus den Färbungen der beiden Blöcke erhalten wir Block A1 (nächste Seite):

| | $E_{g_1}$ / $E_{r_1}$ | $E_{r_2}$ | $E_{r_3}$ | $E_{g_3}$ / $E_{r_4}$ | $E_{g_2}$ |
|---|---|---|---|---|---|
| | | | | 1 | |
| | | | | | 2 |
| | | | | 4 | |
| G | 1 | | | | |
| | | | | 1 | |
| | | | | | 4 |
| | | | | 2 | |
| | | | 1 | 4 | |
| | | 2 | | | |
| | | | 4 | 1 | |
| R, B1 | 1 | | | | 2 |
| | | | 1 | 4 | |
| | | 4 | | | |
| | | | 2 | 1, 4 | |
| | | | 1 | 2 | |
| | | 4 | | | |
| | | | 2 | 1 | |
| R, B2 | 1 | | | | 4 |
| | | | 1 | 2 | |
| | | 2 | | | |
| | | | 4 | 1, 2 | |

| $E_{r_1}^{s_1}$ | $E_{r_2}$ | $E_{r_3}$ | $E_{r_4}^{s_3}$ |
|---|---|---|---|
| | | | 2 |
| | | 1 | |
| | | | 4 |
| | 2 | | |
| | | | 2 |
| | | 4 | |
| | | | 1 |
| G+R, A1 | | | 4 |
| | | 1 | |
| | | | 2 |
| | 4 | | |
| | | | 4 |
| | | 2 | |
| | | | 1 |

Block B6. Aus den Färbungen dieser beiden Blöcke erhalten wir Block A3:

| | $E_{g_1}$ | | | $E_{g_3}$ | $E_{g_2}$ |
|---|---|---|---|---|---|
| | $E_{r_1}$ | $E_{r_2}$ | $E_{r_3}$ | $E_{r_4}$ | |
| | | | | 4 | |
| | | | | | 1 |
| G | 4 | | | | |
| | | | | 4 | |
| | | | | | 2 |
| | | | | 1 | |
| | | | 4 | 2 | |
| | | 1 | | | |
| | | | 2 | 4 | |
| R, B6 | 4 | | | | 1 |
| | | | 4 | 2 | |
| | | 2 | | | |
| | | | 1 | 4, 2 | |
| | | | 4 | 1 | |
| | | 1 | | | |
| | | | 2 | 4, 1 | |
| R, B5 | 4 | | | | 2 |
| | | | 4 | 1 | |
| | | 2 | | | |
| | | | 1 | 4 | |
| | | | | 1 | |
| | | | 4 | | |
| | | | | 2 | |
| | | 1 | | | |
| | | | | 1 | |
| | | | 2 | | |
| | | | | 4 | |
| G+R | 4 | | | | |
| | | | | 2 | |
| | | | 4 | | |
| | | | | 1 | |
| | | 2 | | | |
| | | | | 2 | |
| | | | 1 | | |
| | | | | 4 | |

Block B4. Aus den Färbungen dieser beiden Blöcke erhalten wird Block A2.

| | $E_{g_1}$ / $E_{r_1}$ | $E_{r_2}$ | $E_{r_3}$ | $E_{g_3}$ / $E_{r_4}$ | $E_{g_2}$ |
|---|---|---|---|---|---|
| | | | | 2 | |
| | | | | | 1 |
| | | | | 4 | |
| G | 2 | | | | |
| | | | | 2 | |
| | | | | | 4 |
| | | | | 1 | |
| | | | 2 | 4 | |
| | | 1 | | | |
| | | | 4 | 2 | |
| B3 | 2 | | | | 1 |
| | | | 2 | 4 | |
| | | 4 | | | |
| | | | 1 | 2, 4 | |
| | | | 2 | 1 | |
| | | 4 | | | |
| | | | 1 | 2 | |
| B4 | 2 | | | | 4 |
| | | | 2 | 1 | |
| | | 1 | | | |
| | | | 4 | 2, 1 | |
| | | | | 1 | |
| | | | | 4 | |
| | | 1 | | | |
| | | | | 1 | |
| | | | 4 | | |
| | | | | 2 | |
| G+R, A2 | 2 | | | | |
| | | | | 4 | |
| | | | 2 | | |
| | | | | 1 | |
| | | 4 | | | |
| | | | | 4 | |
| | | | 1 | | |
| | | | | 2 | |

Für $E_{g_1} = E_2$ ist die Konstruktion im wesentlichen dieselbe, für $E_{g_1}$ ist dann nur $f(E_{g_1}) = 2$ zu betrachten.

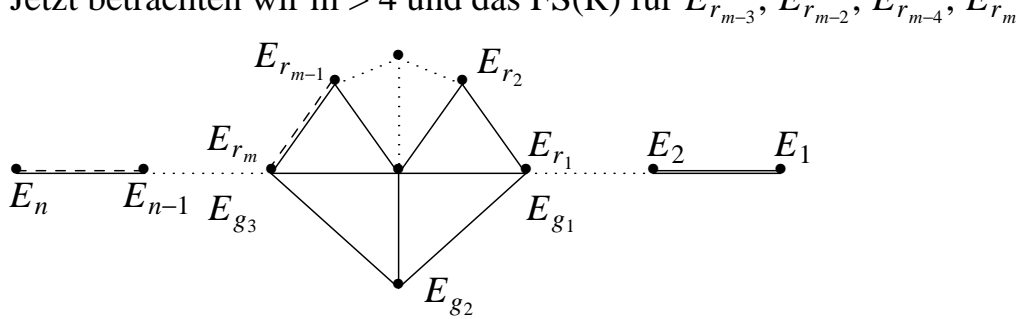

Damit wir ein FS für G+R erhalten, muß - wie im Fall m=4 - bei der Konstruktion der Block des FS(R), der die Kopplung $E_{r_{m-1}}, E_{r_m}$ beschreibt, durch einen Block ohne Kopplung ersetzt werden. Welche Blöcke im FS(R) tatsächlich vorliegen hängt von m ab.

Für $f(E_{g_1}) = 2$, $f(E_{g_2}) = 1$, und die dazu möglichen Blöcken B3, B5, C1 im FS(R) (s. oben) erhalten wir als jeweils möglichen Block im FS(R) zu $f(E_{g_2}) = 4$
Block B3 => $f(E_{r_{m-3}}) = 2$ => zweiter Block ist B4, erledigt s. oben
Block B5 => $f(E_{r_{m-3}}) = 4$ => zweiter Block ist C5
Block C1 => $f(E_{r_{m-3}}) = 1$ => zweiter Block ist B2

Für $f(E_{g_2}) = 1$ ist $f(E_{g_3}) = 2$ und 4, für $f(E_{g_2}) = 4$ ist $f(E_{g_3}) = 1$ und 2. Dies liefert die Bedingungen dafür, welche Färbungen für G+R möglich sind wegen $E_{r_m} = E_{g_3}$.
Entsprechendes gilt für die nächsten Fälle.

| | $E_{r_{m-3}}$ | $E_{r_{m-2}}$ | $E_{r_{m-1}}$ | $E_{r_m}$ | $E_{r_{m-3}}$ | $E_{r_{m-2}}$ | $E_{r_{m-1}}$ | $E_{r_m}$ | |
|---|---|---|---|---|---|---|---|---|---|
| | | | 4 | 2 | | | | | 1 | |
| | | 1 | | | | | | 4 | 2 | |
| R, B5 | 4 | | | 2 | | 1 | | | 1 | |
| | | | 4 | 2 | | | | | 1 | |
| | | 2 | | | | | 2 | | | |
| | | 1 | 4, 2 | | 4 | | | | | G+R, A3 |
| | | | 4 | 2, 1 | | | | | 2 | |
| | | 1 | | | | | | 4 | 1 | |
| R, C5 | 4 | | | 2 | | 2 | | | 1 | |
| | | | 4 | 1, 2 | | | | | 2 | |
| | | 2 | | | | | 1 | | | |
| | | 1 | 2 | | | | | | 4 | |

| $E_{r_{m-3}}$ | $E_{r_{m-2}}$ | $E_{r_{m-4}}$ | $E_{r_m}$ | $E_{r_{m-3}}$ | $E_{r_{m-2}}$ | $E_{r_{m-1}}$ | $E_{r_m}$ |
|---|---|---|---|---|---|---|---|
| | | 1 | 4, 2 | | | | 2 |
| | 2 | | | | | 1 | |
| | | 4 | 2 | | | | 4 |
| R, C1  1 | | | | | 2 | | |
| | | 1 | 2, 4 | | | | 2 |
| | 4 | | | | | 4 | |
| | | 2 | 4 | | | | 1 |
| | | | | 1 | | | G+R, A1 |
| | | 1 | 2 | | | | 4 |
| | 4 | | | | | 1 | |
| | | 2 | 1 | | | | 2 |
| R, B2  1 | | | | | 4 | | |
| | | 1 | 2 | | | | 4 |
| | 2 | | | | | 2 | |
| | | 4 | 1, 2 | | | | 1 |

Für $f(E_{g_1}) = 1$, $f(E_{g_2}) = 2$, und die dazu möglichen Blöcken (s. oben) erhalten wir als jeweils möglichen Block (im FS(R) zu $f(E_{g_2}) = 4$)

Block B1 => $f(E_{r_{m-3}}) = 2$ => zweiter Block ist B2, erledigt s. oben
Block B6 => $f(E_{r_{m-3}}) = 4$ => zweiter Block ist C5
Block C3 => $f(E_{r_{m-3}}) = 1$ => zweiter Block ist B4

| $E_{r_{m-3}}$ | $E_{r_{m-2}}$ | $E_{r_{m-4}}$ | $E_{r_m}$ | $E_{r_{m-3}}$ | $E_{r_{m-2}}$ | $E_{r_{m-1}}$ | $E_{r_m}$ |
|---|---|---|---|---|---|---|---|
| | | 4 | 1 | | | | 1 |
| | 2 | | | | | 4 | |
| | | 1 | 4 | | | | 2 |
| R, B6  4 | | | | | 1 | | |
| | | 4 | 1 | | | | 1 |
| | 1 | | | | | 2 | |
| | | 2 | 4, 1 | | | | 4 |
| | | | | 4 | | | |
| | | 4 | 2, 1 | | | | 2 |
| | 1 | | | | | 4 | |
| | | 2 | 1 | | | | 1 |
| R, C5  4 | | | | | 2 | | |
| | | 4 | 1, 2 | | | | 2 |
| | 2 | | | | | 1 | |
| | | 1 | 2 | | | | 4 |

|  |  | $r_{m-3}$ | $r_{m-2}$ | $r_{m-4}$ | $r_m$ | $r_{m-3}$ | $r_{m-2}$ | $r_{m-1}$ | $r_m$ |
|---|---|---|---|---|---|---|---|---|---|
|  |  |  |  | 2 | 4, 1 |  |  |  | 1 |
|  |  |  | 1 |  |  |  |  | 2 |  |
|  |  |  |  | 4 | 1 |  |  |  | 4 |
| C3 | 2 |  |  |  |  |  | 1 |  |  |
|  |  |  |  | 2 | 1, 4 |  |  |  | 1 |
|  |  |  | 4 |  |  |  |  | 4 |  |
|  |  |  |  | 1 | 4 |  |  |  | 2 |
|  |  |  |  |  |  | 2 |  |  |  |
|  |  |  |  | 2 | 1 |  |  |  | 4 |
|  |  |  | 4 |  |  |  |  | 2 |  |
|  |  |  |  | 1 | 2 |  |  |  | 1 |
| B4 | 2 |  |  |  |  | 4 |  |  |  |
|  |  |  |  | 2 | 1 |  |  |  | 4 |
|  |  |  | 1 |  |  |  |  | 1 |  |
|  |  |  |  | 4 | 2, 1 |  |  |  | 2 |

Für $f(E_{g_1}) = 4$, $f(E_{g_2}) = 1$, und die dazu möglichen Blöcken (s. oben) erhalten wir als jeweils möglichen Block (im FS(R) zu $f(E_{g_2}) = 2$)

Block B5 => $f(E_{r_{m-3}}) = 4$ => zweiter Block ist B6, erledigt s. oben

Block B3 => $f(E_{r_{m-3}}) = 2$ => zweiter Block ist C3

Block C1 => $f(E_{r_{m-3}}) = 1$ => zweiter Block ist B1

|  |  | $E_{r_{m-3}}$ | $E_{r_{m-2}}$ | $E_{r_{m-4}}$ | $E_{r_m}$ | $E_{r_{m-3}}$ | $E_{r_{m-2}}$ | $E_{r_{m-1}}$ | $E_{r_m}$ |
|---|---|---|---|---|---|---|---|---|---|
|  |  |  |  | 2 | 4 |  |  |  | 1 |
|  |  |  | 1 |  |  |  |  | 2 |  |
|  |  |  |  | 4 | 2 |  |  |  | 4 |
| R, B3 | 2 |  |  |  |  |  | 1 |  |  |
|  |  |  |  | 2 | 4 |  |  |  | 1 |
|  |  |  | 4 |  |  |  |  | 4 |  |
|  |  |  |  | 1 | 2, 4 |  |  |  | 2 |
|  |  |  |  |  |  | 2 |  |  |  |
|  |  |  |  | 2 | 4, 1 |  |  |  | 4 |
|  |  |  | 1 |  |  |  |  | 2 |  |
|  |  |  |  | 4 | 1 |  |  |  | 1 |
| R, C3 | 2 |  |  |  |  | 4 |  |  |  |
|  |  |  |  | 2 | 1, 4 |  |  |  | 4 |
|  |  |  | 4 |  |  |  |  | 1 |  |
|  |  |  |  | 1 | 4 |  |  |  | 2 |

| | | $r_{m-3}$ | $r_{m-2}$ | $r_{m-4}$ | $r_m$ | | $r_{m-3}$ | $r_{m-2}$ | $r_{m-1}$ | $r_m$ |
|---|---|---|---|---|---|---|---|---|---|---|
| | | | | 1 | 4, 2 | | | | | 2 |
| | | | 2 | | | | | | 1 | |
| | | | | 4 | 2 | | | | | 4 |
| R, C1 | 1 | | | | | | 2 | | | |
| | | | | 1 | 2, 4 | | | | | 2 |
| | | | 4 | | | | | 4 | | |
| | | | | 2 | 4 | | | | | 1 |
| | | | | | | 1 | | | | |
| | | | | 1 | 4 | | | | | 4 |
| | | | 2 | | | | | | 1 | |
| | | | | 4 | 1 | | | | | 2 |
| R, B1 | 1 | | | | | | 4 | | | |
| | | | | 1 | 4 | | | | | 4 |
| | | | 4 | | | | | 2 | | |
| | | | | 2 | 1, 4 | | | | | 1 |

Die konstruierten Färbungen bilden ein FS für G+R, insbesondere gibt es nur die gekoppelten Spalten der Ecken $E_{n-1}, E_n$.

von G+R. $E_{n-1} = E_{g_{k-1}}$ fällt nach innen.

2.1. k=3, die gekoppelte Ecken $E_{r_m}$, $E_{r_{m-1}}$ von R liegen auf dem Rand von G+R.

2.1.1 m=3, (ein Fall m=k) wir setzen $f(E_{r_2}) = f(E_{g_2})$.

2.1.2 m = 4, $E_{g_1} > E_2$

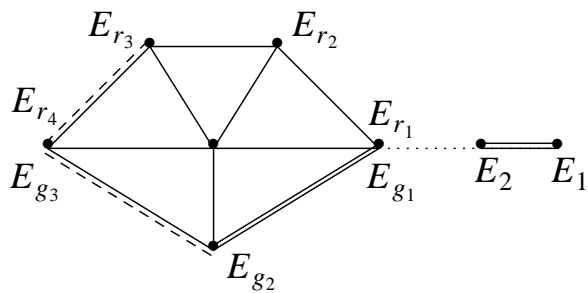

$E_1$, $E_2$ Basisseite von G, $f(E_1) = 1$, $f(E_2) = 2$.  $E_{g_1}$, $E_{g_2}$ Basisseite von R.

$E_{g_2}$, $E_{g_3}$ gekoppelte Spalten von G, $E_{r_3}$, $E_{r_4}$ gekoppelte Spalten von R

Für eine Farbe von $E_{g_1}$ haben wir zwei Farben für $E_{g_2}$ und daher zwei FS für R, da $E_{g_1}$, $E_{g_2}$ Basisseite von R. Wir betrachten alle möglichen Fälle und konstruieren dazu jeweils das FS von G + R.

$f(E_{g_1}) = 1$, $f(E_{g_2}) = 2\, oder\, 4$. Block von FS(G) (nur ein "halber" Block, die andere Hälfte gehört zu einem anderen Wert für $f(E_{g_1})$. Von folgenden Blöcken kommt der erste in Frage wegen $f(E_{g_3}) \neq 1$).

| $E_{g_1}$ | $E_{g_2}$ | $E_{g_3}$ | | $E_{g_1}$ | $E_{g_2}$ | $E_{g_3}$ | | $E_{g_1}$ | $E_{g_2}$ | $E_{g_3}$ |
|---|---|---|---|---|---|---|---|---|---|---|
| | 2 | 4 | | | 2 | 4, 1 | | | 4 | 2, 1 |
| 1 | | | | 1 | | | | 1 | | |
| | 4 | 2 | | | 4 | 1 | | | 2 | 1 |

dazu die Blöcke des FS(R) abhängig von $f(E_{g_1})$, $f(E_{g_2})$. Die Spalte $f(E_{g_3})$ zeigt die Bedingung für passende Färbungen.

| $g_2$ | $g_1$ | $r_2$ | $r_3$ | $r_4$ | $g_3$ |
|---|---|---|---|---|---|
| | | | 1 | 4 | |
| | | 2 | | | |
| | | | 4 | 1 | |
| 2 | 1 | | | | 4 |
| | | | 1 | 4 | |
| | | 4 | | | |
| | | | 2 | 1, 4 | |
| | | | 1 | 2 | |
| | | 2 | | | |
| | | | 4 | 1, 2 | |
| 4 | 1 | | | | 2 |
| | | | 1 | 2 | |
| | | 4 | | | |
| | | | 2 | 1 | |

Die möglichen Färbungen auf dem Rand von G+R (unter der Bedingung $f(E_{g_3}) = f(E_{r_4})$) sind:

| $E_{r_1}$ | $E_{r_2}$ | $E_{r_3}$ | $E_{r_4}$ |
|---|---|---|---|
| 1 | 2 | 1 | 4 |
| 1 | 4 | 1 | 4 |
| 1 | 4 | 2 | 4 |
| 1 | 2 | 1 | 2 |
| 1 | 2 | 4 | 2 |
| 1 | 4 | 1 | 2 |

die möglichen Färbungen für FS(G+R) als Block

| $E_{r_1}$ | $E_{r_2}$ | $E_{r_3}$ | $E_{r_4}$ |
|---|---|---|---|
| | | 1 | 4, 2 |
| | 2 | | |
| | | 4 | 2 |
| 1 | | | |
| | | 1 | 2, 4 |
| | 4 | | |
| | | 2 | 4 |

ein Block C1.

des FS(R) abhängig von $f(E_{g_2})=1$ bzw. $f(E_{g_2})=4$.

| | $E_{g_1}$ / $E_{r_1}$ | $E_{r_2}$ | $E_{r_3}$ | $E_{g_3}$ / $E_{r_4}$ | $E_{g_2}$ / $E_{r_5}$ |
|---|---|---|---|---|---|
| | | | | 4, 2 | 1 |
| G | 2 | | | | |
| | | | | 2 | 4 |
| | | | 2 | 4 | |
| | | 1 | | | |
| | | | 4 | 2 | |
| R, B3 | 2 | | | | 1 |
| | | | 2 | 4 | |
| | | 4 | | | |
| | | | 1 | 2, 4 | |
| | | | 2 | 1 | |
| | | 1 | | | |
| | | | 4 | 2, 1 | |
| R, B4 | 2 | | | | 4 |
| | | | 2 | 1 | |
| | | 4 | | | |
| | | | 1 | 2 | |

Die möglichen Färbungen auf dem Rand von G+R sind:

| | | | |
|---|---|---|---|
| 2 | 1 | 2 | 4 |
| 2 | 1 | 4 | 2 |
| 2 | 4 | 2 | 4 |
| 2 | 4 | 1 | 2, 4 |
| 2 | 1 | 4 | 2 |
| 2 | 4 | 1 | 2 |

die passenden Färbungen als Block B3

| | | | | |
|---|---|---|---|---|
| | | | 2 | 4 |
| | 1 | | | |
| | | | 4 | 2 |
| 2 | | | | |
| | | | 2 | 4 |
| | 4 | | | |
| | | | 1 | 2, 4 |

(Von Block R, B3 passen übrigens schon alle Färbungen.)

$\gamma(E_{g_1})$ ... = 4, $\gamma(E_{g_2})$ = 1, 2. Der Block von FS(G) und dazu die Blöcke B5 bzw. B6 des FS(R) abhängig von $f(E_{g_2})=1$ bzw. $f(E_{g_2})=2$.

| | $E_{g_1}$ | | | $E_{g_3}$ | $E_{g_2}$ |
|---|---|---|---|---|---|
| | $E_{r_1}$ | $E_{r_2}$ | $E_{r_3}$ | $E_{r_4}$ | $E_{r_5}$ |
| | | | | 1 | 2, 4 |
| G | 4 | | | | |
| | | | | 2 | 4 |
| | | | 2 | 4 | |
| | | 1 | | | |
| | | | 4 | 2 | |
| R, B5 | 4 | | | | 2, 4 |
| | | | 1 | 4, 2 | |
| | | 2 | | | |
| | | | 4 | 2 | |
| | | | 4 | 1 | |
| | | 1 | | | |
| | | | 2 | 4, 1 | |
| R, B6 | 4 | | | | 4 |
| | | | 4 | 1 | |
| | | 2 | | | |
| | | | 1 | 4 | |

Die möglichen Färbungen auf dem Rand von G+R sind:

| | | | |
|---|---|---|---|
| 4 | 1 | 2 | 4 |
| 4 | 1 | 4 | 2 |
| 4 | 2 | 1 | 4, 2 |
| 4 | 2 | 4 | 2 |
| 4 | 1 | 2 | 4 |
| 4 | 2 | 1 | 4 |

die passenden Färbungen als Block B3

| | | | |
|---|---|---|---|
| | | 4 | 2 |
| | 1 | | |
| | | 2 | 4 |
| 4 | | | |
| | | 4 | 2 |
| | 2 | | |
| | | 1 | 4, 2 |

Die Anzahl der Färbungen nimmt zu entsprechend $\gamma(G + R) = \gamma(G) + 1$. Die gekoppelten Spalten sind die Spalten zu $E_{r_3}$, $E_{r_4} = E_n$.

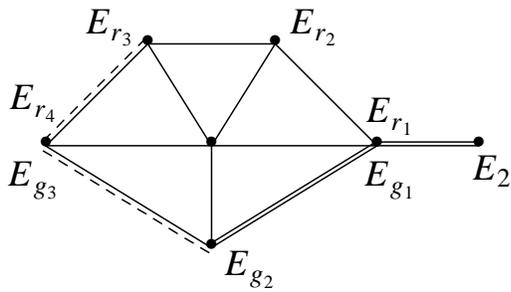

$E_{g_1} = E_2 \Rightarrow f(E_{g_1}) = f(E_2) = 2.$

Wie 2.1.2 mit $f(E_{g_1}) = 2$.

|  | $E_1$ | $E_{g_1}$ / $E_{r_1}$ | $E_{r_2}$ | $E_{r_3}$ | $E_{g_3}$ / $E_{r_4}$ | $E_{g_2}$ |
|---|---|---|---|---|---|---|
| G | 1 | 2 |  |  | 4, 2 | 1 |
|  |  |  |  |  | 2 | 4 |
|  |  |  |  | 2 | 4 |  |
|  |  |  | 1 |  |  |  |
|  |  |  |  | 4 | 2 |  |
| R, B3 |  | 2 |  |  |  | 1 |
|  |  |  |  | 2 | 4 |  |
|  |  |  | 4 |  |  |  |
|  |  |  |  | 1 | 2, 4 |  |
|  |  |  |  | 2 | 1 |  |
|  |  |  | 4 |  |  |  |
|  |  |  |  | 1 | 2 |  |
| R, B4 |  | 4 |  |  |  | 4 |
|  |  |  |  | 4 | 1 |  |
|  |  |  | 2 |  |  |  |
|  |  |  |  | 1 | 4 |  |

die passenden Färbungen als Block B3

|  | $E_1$ | $E_{g_1}$ / $E_{r_1}$ | $E_{r_2}$ | $E_{r_3}$ | $E_{g_3}$ / $E_{r_4}$ | $E_{g_2}$ |
|---|---|---|---|---|---|---|
|  |  |  |  | 2 | 4 |  |
|  |  |  | 1 |  |  |  |
|  | 1 | 2 |  |  |  |  |
|  |  |  |  | 4 | 2 |  |
|  |  |  |  | 2 | 4 |  |
|  |  |  | 4 |  |  |  |
|  |  |  |  | 1 | 2, 4 |  |

Die Färbungen das FS(R), B4 mit $f(E_{r_4}) = 1$ fallen weg, da $E_{r_4}$ ( $= E_{g_3}$) mit $E_1$ durch Kante verbunden ist und $f(E_1) = 1$.

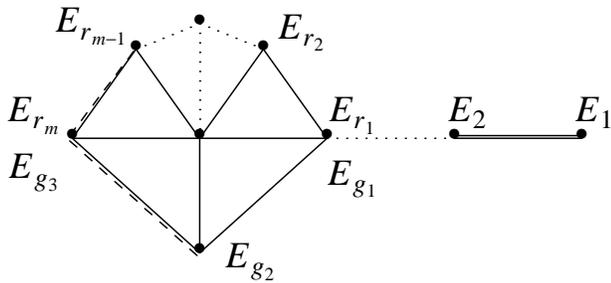

$E_1, E_2$ Basisseite von G, $f(E_1) = 1$, $f(E_2) = 2$. $E_{g_1}, E_{g_2}$ Basisseite von R.

$E_{g_2}, E_{g_3}$ gekoppelte Spalten von G, $E_{r_{m-1}}, E_{r_m}$ gekoppelte Spalten von R

$f(E_{g_1}) = 1 \Rightarrow f(E_{g_2}) = 2, 4$

$f(E_{g_1}) = 2 \Rightarrow f(E_{g_2}) = 1, 4$

$f(E_{g_1}) = 4 \Rightarrow f(E_{g_2}) = 1, 2$

Zunächst sind für $E_{r_{m-3}}$, $E_{r_{m-2}}$, $E_{r_{m-1}}$, $E_{r_m}$ im FS(R) folgende Blöcke FS(R) möglich:

$f(E_{g_2}) = 1$ : B3, B5, C1

$f(E_{g_2}) = 2$ : B1, B6, C3

$f(E_{g_2}) = 4$ : B2, B4, C5

$f(E_{g_3})$ hat nur die Werte 2, 4. $(f(E_{g_3}) \neq f(E_1) = 1)$

Zu zwei Werten von $f(E_{g_2})$ gibt es je ein FS(R) beginnend mit der Basisseite $E_{g_2} E_{g_1}$, die gleich sind außer in der Spalte, die von $E_{g_2}$ abhängt. Haben wir $f(E_{g_2}) = 2$ und Block B1, so gibt es dasselbe FS bis auf den letzten Block, dieser muß B2 sein. (Die ersten drei Spalten in B1, B2 müssen gleich sein.) So ergeben sich die folgenden möglichen Kombinationen:

| $f(E_{g_1}) = 1$ | | $f(E_{g_1}) = 2$ | | $f(E_{g_1}) = 4$ | |
|---|---|---|---|---|---|
| $f(E_{g_2}) = 2, 4$ | | $f(E_{g_2}) = 1, 4$ | | $f(E_{g_2}) = 1, 2$ | |
| B1 | mit B2 | B3 | mit B4 | B3 | mit C3 |
| B6 | C5 | B5 | C5 | B5 | B6 |
| C3 | B4 | C1 | B2 | C1 | B1 |

Für diese Kombinationen konstruieren wir den Block zu $E_{r_{m-3}}$, $E_{r_{m-2}}$, $E_{r_{m-1}}$, $E_{r_m}$ im FS(G+R).

$f(E_{g_1}) = 1, f(E_{g_2}) = 2, 4$

| | $E_{r_{m-3}}$ | $E_{r_{m-2}}$ | $E_{r_{m-1}}$ | $E_{g_3}$ / $E_{r_m}$ | $E_{g_2}$ | $E_{g_3}$ |
|---|---|---|---|---|---|---|
| | | | | 4, 2 | 2 | |
| G | | | | | | 1 |
| | | | | 2 | 4 | |
| | | | 1 | 4 | | |
| | | 2 | | | | |
| | | | 4 | 1 | | |
| R, B1 | 1 | | | | 2 | |
| | | | 1 | 4 | | |
| | | 4 | | | | |
| | | | 2 | 1, 4 | | |
| | | | 1 | 2 | | |
| | | 2 | | | | |
| | | | 4 | 1, 2 | | |
| R, B2 | 1 | | | | 4 | |
| | | | 1 | 2 | | |
| | | 4 | | | | |
| | | | 2 | 1 | | |
| passende Färbungen | | | | | | |
| | 1 | 2 | 1 | 4 | | |
| | 1 | 4 | 1 | 4 | | |
| | 1 | 4 | 2 | 4 | | |
| | 1 | 2 | 1 | 2 | | |
| | 1 | 2 | 4 | 2 | | |
| | 1 | 3 | 1 | 2 | | |
| als Block | | | | | | |
| | | | 1 | 4, 2 | | |
| | | 2 | | | | |
| | | | 4 | 2 | | |
| C1 | 1 | | | | | |
| | | | 1 | 4, 2 | | |
| | | 4 | | | | |
| | | | 2 | 4 | | |

$f(E_{g_1}) = 2, f(E_{g_2}) = 1, 4$

| | $E_{r_{m-3}}$ | $E_{r_{m-2}}$ | $E_{r_{m-1}}$ | $E_{r_m}$ / $E_{g_3}$ | $E_{g_2}$ | $E_{g_1}$ |
|---|---|---|---|---|---|---|
| | | | | 4 | 2 | |
| G | | | | | | 1 |
| | | | | 2 | 4 | |
| | | | 2 | 4, 1 | | |
| | | 1 | | | | |
| | | | 4 | 1 | | |
| R, C3 | 2 | | | | 2 | |
| | | | 2 | 1, 4 | | |
| | | 4 | | | | |
| | | | 1 | 4 | | |
| | | | 2 | 1 | | |
| | | 4 | | | | |
| | | | 1 | 2 | | |
| R, B4 | 2 | | | | 4 | |
| | | | 2 | 1 | | |
| | | 1 | | | | |
| | | | 4 | 2, 1 | | |
| passende Färbungen | | | | | | |
| | 2 | 1 | 2 | 4 | | |
| | 2 | 4 | 2 | 4 | | |
| | 2 | 4 | 1 | 4 | | |
| | 2 | 4 | 1 | 2 | | |
| | 2 | 1 | 4 | 2 | | |
| als Block | | | | | | |
| | | | 2 | 4 | | |
| | | 1 | | | | |
| | | | 4 | 2 | | |
| B3 | 2 | | | | | |
| | | | 2 | 4 | | |
| | | 4 | | | | |
| | | | 1 | 4, 2 | | |

$f(E_{g_1}) = 4$, $f(E_{g_2}) = 1, 2$

| | $E_{r_{m-3}}$ | $E_{r_{m-2}}$ | $E_{r_{m-1}}$ | $E_{r_m}$ / $E_{g_3}$ | $E_{g_2}$ | $E_{g_1}$ |
|---|---|---|---|---|---|---|
| | | | | 2, 4 | 1 | |
| G | | | | | | 4 |
| | | | | 4 | 2 | |
| | | | 4 | 1 | | |
| | | 2 | | | | |
| | | | 1 | 4 | | |
| R, B6 | 4 | | | | 2 | |
| | | | 4 | 1 | | |
| | | 1 | | | | |
| | | | 2 | 4, 1 | | |
| | | | 4 | 2, 1 | | |
| | | 1 | | | | |
| | | | 2 | 1 | | |
| R, C5 | 4 | | | | 4 | |
| | | | 4 | 1, 2 | | |
| | | 2 | | | | |
| | | | 1 | 2 | | |
| passende Färbungen | | | | | | |
| | 4 | 2 | 1 | 4 | | |
| | 4 | 1 | 2 | 4 | | |
| | 4 | 1 | 4 | 2 | | |
| | 4 | 2 | 4 | 2 | | |
| | 4 | 2 | 1 | 2 | | |
| als Block | | | | | | |
| | | | 2 | 4 | | |
| | | 1 | | | | |
| | | | 4 | 2 | | |
| B5 | 4 | | | | | |
| | | | 1 | 4, 2 | | |
| | | 2 | | | | |
| | | | 4 | 2 | | |

Damit ist für $f(E_{g_1}) = 1$, $f(E_{g_2}) = 2, 4$ für die drei möglichen Fälle B1 mit B2, B6 mit C5 und C3 mit B4 das FS(G+R) konstruiert.

Nun zu den jeweils drei möglichen Fällen für $f(E_{g_1}) = 2$, $f(E_{g_2}) = 1, 4$ und $f(E_{g_1}) = 4$, $f(E_{g_2}) = 1, 2$.

$f(E_{g_2}) = 1$ bzw. 4, $f(E_{g_3}) = 2, 4$ bzw. 2

| | $E_{g_1}$ | $E_{g_2}$ | $E_{g_3}$ |
|---|---|---|---|
| | | 1 | 2, 4 |
| G | 2 | | |
| | 4 | 2 | |

| | $E_{r_{m-3}}$ | $E_{r_{m-2}}$ | $E_{r_{m-1}}$ | $E_{r_m}$ | $E_{g_3}$ |
|---|---|---|---|---|---|
| | | | 2 | 4 | |
| | | 1 | | | |
| | | | 4 | 2 | |
| B3 | 2 | | | | 2, 4 |
| | | | 2 | 4 | |
| | | 4 | | | |
| | | | 1 | 2, 4 | |

| | $E_{r_{m-3}}$ | $E_{r_{m-2}}$ | $E_{r_{m-1}}$ | $E_{r_m}$ | $E_{g_3}$ |
|---|---|---|---|---|---|
| | | | 2 | 1 | |
| | | 4 | | | |
| | | | 1 | 2 | |
| B4 | 2 | | | | 2 |
| | | | 2 | 1 | |
| | | 1 | | | |
| | | | 4 | 2, 1 | |

Der Block B3 kann komplett für FS(G+R) genommen werden.

Von den Blöcken B5, C5 kann B5 komplett für FS(G+R) genommen werden:

| | $E_{r_{m-3}}$ | $E_{r_{m-2}}$ | $E_{r_{m-1}}$ | $E_{r_m}$ | $E_{g_2}$ |
|---|---|---|---|---|---|
| | | | 4 | 2 | |
| | | 1 | | | |
| | | | 2 | 4 | |
| B5 | 4 | | | | 2, 4 |
| | | | 4 | 2 | |
| | | 2 | | | |
| | | | 1 | 4, 2 | |

Von den Blöcken C1, B2 kann C1 komplett für FS(G+R) genommen werden:

| | $E_{r_{m-3}}$ | $E_{r_{m-2}}$ | $E_{r_{m-1}}$ | $E_{r_m}$ | $E_{g_2}$ |
|---|---|---|---|---|---|
| | | | 1 | 4, 2 | |
| | | 2 | | | |
| | | | 4 | 2 | |
| C1 | 1 | | | | 2, 4 |
| | | | 1 | 2, 4 | |
| | | 4 | | | |
| | | | 2 | 4 | |

$f(E_{g_2}) = 1$ bzw. 2, $f(E_{g_3}) = 2, 4$ bzw. 4

| | $E_{g_2}$ | $E_{g_3}$ |
|---|---|---|
| | 1 | 2, 4 |
| G | 4 | |
| | 2 | 4 |

B3, B5, C1 können komplett für FS(G+R) genommen werden.

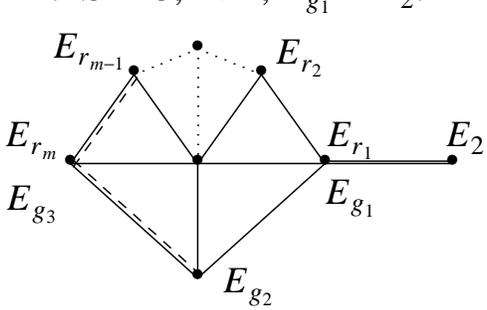

Wir erhalten wie in 2.1.4 mit $f(E_{g_1}) = 2$ für FS(G+R) den Block B3.

Weitere Färbungen von $f(E_{g_2}) = 4$ sind möglich, werden nicht benötigt.

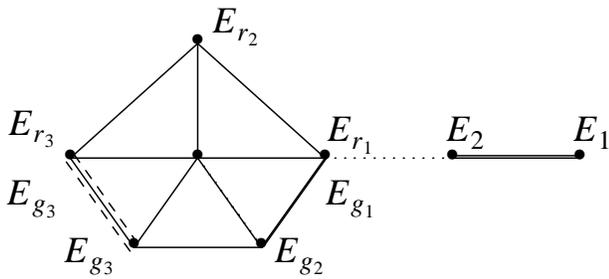

$E_{g_4}$, $E_{g_3}$ gekoppelte Ecken von G. $E_{r_3}$, $E_{r_4}$ (= $E_{g_3}$) gekoppelte Ecken von R. Die gekoppelten Ecken von G werden identifiziert mit den gekoppelten Ecken von R.

$f(E_{g_1}) = 1, 2$ oder $4$.
$f(E_{g_1}) = 1 => f(E_{g_2}) = 2, 4$.

Block FS(G)

|       | $E_{g_1}$ | $E_{g_2}$ | $E_{g_3}$ | $E_{g_4}$ |
|-------|-----------|-----------|-----------|-----------|
|       |           |           | 1         | 4, 2      |
|       |           | 2         |           |           |
|       |           |           | 4         | 2         |
| G, C1 | 1         |           |           |           |
|       |           |           | 1         | 2, 4      |
|       |           | 4         |           |           |
|       |           |           | 2         | 4         |

Block FS(R) für $f(E_{g_2}) = 2$ und $f(E_{g_2}) = 4$.

|       | $E_{r_5}$ | $E_{r_1}$ | $E_{r_2}$ | $E_{r_3}$ | $E_{r_4}$ |
|-------|-----------|-----------|-----------|-----------|-----------|
|       |           |           |           | 1         | 4         |
|       |           |           | 2         |           |           |
|       |           |           |           | 4         | 1         |
| R, B1 | 2         | 1         |           |           |           |
|       |           |           |           | 1         | 4         |
|       |           |           | 4         |           |           |
|       |           |           |           | 2         | 1, 4      |

|       | $E_{r_5}$ | $E_{r_1}$ | $E_{r_2}$ | $E_{r_3}$ | $E_{r_4}$ |
|-------|-----------|-----------|-----------|-----------|-----------|
|       |           |           |           | 1         | 2         |
|       |           |           | 2         |           |           |
|       |           |           |           | 4         | 1, 2      |
| R, B2 | 4         | 1         |           |           |           |
|       |           |           |           | 1         | 2         |
|       |           |           | 4         |           |           |
|       |           |           |           | 2         | 1         |

in einer Tabelle

| $E_{r_1}$ | $E_{r_2}$ | $E_{r_3}$ | $E_{r_4}$ | $E_{r_5}$ |
|---|---|---|---|---|
| $E_{g_1}$ | | $E_{g_4}$ | $E_{g_3}$ | $E_{g_2}$ |
| | | 4, 2 | 1 | |
| | | | | 2 |
| | | 2 | 4 | |
| G 1 | | | | |
| | | 2, 4 | 1 | |
| | | | | 4 |
| | | 4 | 2 | |

| | | 1 | 4 | |
|---|---|---|---|---|
| | 2 | | | |
| | | 4 | 1 | |
| R, B1 1 | | | | 2 |
| | | 1 | 4 | |
| | 4 | | | |
| | | 2 | 1, 4 | |

| | | 1 | 2 | |
|---|---|---|---|---|
| | 2 | | | |
| | | 4 | 1, 2 | |
| R, B2 1 | | | | 4 |
| | | 1 | 2 | |
| | 4 | | | |
| | | 2 | 1 | |

finden wir als passende Färbungen

| 1 | 2 | 4 | 1 | 2 |
|---|---|---|---|---|
| | 4 | 2 | 1 | 2 |
| | 4 | 2 | 4 | 2 |
| | 4 | 2 | 1 | 4 |
| | 2 | 4 | 1 | 4 |
| | 2 | 4 | 2 | 4 |

als Block des FS(G+R)

| | 2 | 4 | | |
|---|---|---|---|---|
| 1 | | | | |
| | 4 | 2 | | |

Zu einem Wert $f(E_{g_1}) = 1$ gibt es die gekoppelten Spalten für $E_{r_2}$ $E_{r_3}$.

Dasselbe für $f(E_{g_1}) = 2$, $f(E_{g_2}) = 1, 4$.

| | $E_{r_1}$ $E_{g_1}$ | $E_{r_2}$ | $E_{r_3}$ $E_{g_4}$ | $E_{r_4}$ $E_{g_3}$ | $E_{r_5}$ $E_{g_2}$ |
|---|---|---|---|---|---|
| | | | 4 | 1 | |
| | | | | | 1 |
| | | | 2 | 4 | |
| G, B3 | 2 | | | | |
| | | | 4 | 2 | |
| | | | | | 4 |
| | | | 2, 4 | 1 | |
| | | | 2 | 4 | |
| | | 1 | | | |
| | | | 4 | 2 | |
| R, B3 | 2 | | | | 1 |
| | | | 2 | 4 | |
| | | 4 | | | |
| | | | 1 | 4 | |
| | | | 2 | 1 | |
| | | 1 | | | |
| | | | 4 | 2, 1 | |
| R, B4 | 2 | | | | 4 |
| | | | 2 | 1 | |
| | | 4 | | | |
| | | | 1 | 2 | |

finden wir als passende Färbungen

| | | | | |
|---|---|---|---|---|
| 2 | 1 | 4 | 2 | 1 |
| | 4 | 2 | 4 | 1 |
| | 1 | 2 | 4 | 1 |
| | 1 | 4 | 2 | 4 |
| | 1 | 2 | 1 | 4 |
| | 4 | 2 | 1 | 4 |
| | 1 | 4 | 1 | 4 |

als Block des FS(G+R)

| | | | |
|---|---|---|---|
| | | 1 | 2, 4 |
| 2 | | | |
| | | 4 | 2 |

|        | $E_{r_1}$ | $E_{r_2}$ | $E_{r_3}$ | $E_{r_4}$ | $E_{r_5}$ |
|--------|-----------|-----------|-----------|-----------|-----------|
|        | $E_{g_1}$ |           | $E_{g_4}$ | $E_{g_3}$ | $E_{g_2}$ |
|        |           |           | 2         | 4         |           |
|        |           |           |           |           | 1         |
|        |           |           | 4         | 2         |           |
| G, B5  | 4         |           |           |           |           |
|        |           |           | 2         | 4         |           |
|        |           |           |           |           | 2         |
|        |           |           | 4, 2      | 1         |           |

|        |           |           | 4         | 2         |           |
|--------|-----------|-----------|-----------|-----------|-----------|
|        |           | 1         |           |           |           |
|        |           |           | 2         | 4         |           |
| R, B5  | 4         |           |           |           | 1         |
|        |           |           | 4         | 2         |           |
|        |           | 2         |           |           |           |
|        |           |           | 1         | 4, 2      |           |

|        |           |           | 4         | 1         |           |
|--------|-----------|-----------|-----------|-----------|-----------|
|        |           | 1         |           |           |           |
|        |           |           | 2         | 4, 1      |           |
| R, B6  | 4         |           |           |           | 2         |
|        |           |           | 4         | 1         |           |
|        |           | 2         |           |           |           |
|        |           |           | 1         | 4         |           |

finden wir als passende Färbungen

|   |   |   |   |   |
|---|---|---|---|---|
| 4 | 1 | 2 | 4 | 1 |
| 2 | 4 | 2 | 1 |   |
| 1 | 2 | 4 | 2 |   |
| 1 | 4 | 1 | 2 |   |
| 2 | 4 | 1 | 2 |   |
| 1 | 2 | 1 | 2 |   |

als Block des FS(G+R)

|   |   |   |      |
|---|---|---|------|
|   |   | 1 | 4, 2 |
| 4 |   |   |      |
|   |   | 2 | 4    |

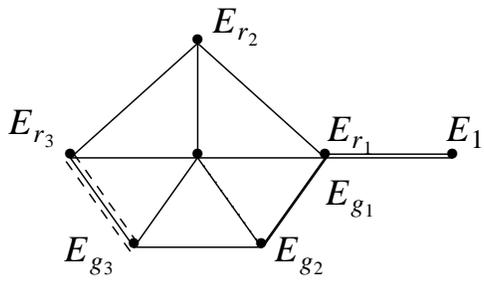

$f(E_{g_1}) = f(E_2) = 2.$
Wie oben, jedoch ist nur der Fall $f(E_{g_1}) = 2$ zu betrachten.

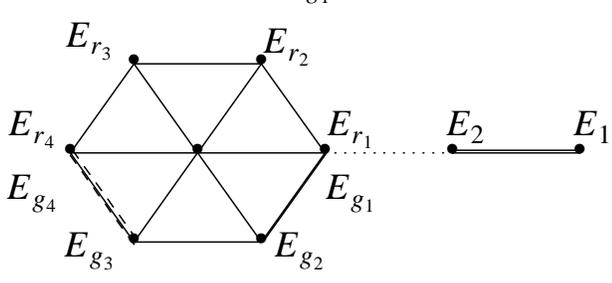

$f(E_{g_1}) = 1.$

| | $E_{g_1}$ | | | $E_{g_4}$ | $E_{g_3}$ | $E_{g_2}$ |
| | $E_{r_1}$ | $E_{r_2}$ | $E_{r_3}$ | $E_{r_4}$ | $E_{r_5}$ | |
|---|---|---|---|---|---|---|
| | | | | 4, 2 | 1 | |
| | | | | | | 2 |
| | | | | 2 | 4 | |
| G | 1 | | | 2, 4 | 1 | |
| | | | | 4 | 2 | 4 |
| | | | | 2 | 4, 1 | |
| | | | 1 | 4 | 1 | |
| | | 2 | | 2 | 1, 4 | |
| | | | 4 | 1 | 4 | |
| R | 1 | | | 4 | 1 | 2 |
| | | | 1 | 2 | 4, 1 | |
| | | 4 | | 4 | 1 | |
| | | | 2 | 1 | 4 | |

Fortsetzung der Tabelle nächste Seite.

| | $E_{g_1}$ | | | $E_{g_4}$ | $E_{g_3}$ | $E_{g_2}$ |
|---|---|---|---|---|---|---|
| | $E_{r_1}$ | $E_{r_2}$ | $E_{r_3}$ | $E_{r_4}$ | $E_{r_5}$ | |
| | | | | 2 | 1 | |
| | | | 1 | | | |
| | | | | 4 | 2, 1 | |
| | | 2 | | | | |
| | | | | 2 | 1 | |
| | | | 4 | | | |
| | | | | 1 | 2 | |
| R | 1 | | | | | 4 |
| | | | | 4 | 2, 1 | |
| | | | 1 | | | |
| | | | | 2 | 1 | |
| | | | 4 | | | |
| | | | | 4 | 1, 2 | |
| | | | 2 | | | |
| | | | | 1 | 2 | |
| | | | 1 | 2, 4 | | |
| | | 2 | | | | |
| | | | 4 | 2 | | |
| G+R, C1 | 1 | | | | | |
| | | | 1 | 4, 2 | | |
| | | | 4 | | | |
| | | | 2 | 4 | | |

k=4, m=3, $E_{g_1} = E_2$

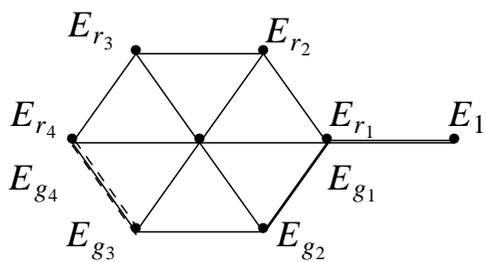

$f(E_{g_1}) = f(E_2) = 2.$
Wie oben, jedoch ist nur der Fall $f(E_{g_1}) = 2$ zu betrachten.

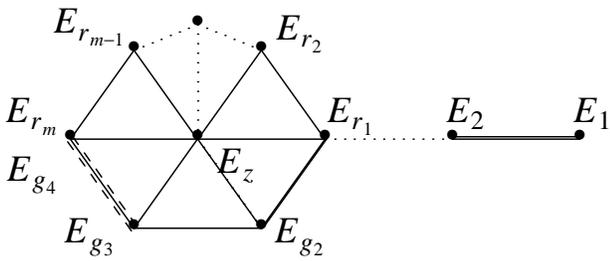

$E_{g_4}$, $E_{g_3}$ gekoppelte Ecken von G. $E_{r_3}$, $E_{r_4}$ (= $E_{g_3}$) gekoppelte Ecken von R.
Die gekoppelten Ecken von G werden identifiziert mit den gekoppelten Ecken von R.

$f(E_{g_1}) = 1$, 2 oder 4.
$f(E_{g_1}) = 1 \Rightarrow f(E_{g_2}) = 2, 4 \Rightarrow$ Block FS(G) = C1.

|        | $E_{g_1}$ | $E_{g_2}$ | $E_{g_3}$ | $E_{g_4}$ |
|--------|-----------|-----------|-----------|-----------|
|        |           |           | 1         | 4, 2      |
|        |           | 2         |           |           |
|        |           |           | 4         | 2         |
| G, C1  | 1         |           |           |           |
|        |           |           | 1         | 2, 4      |
|        |           | 4         |           |           |
|        |           |           | 2         | 4         |

Block FS(R) für $f(E_{g_2}) = 2$ und $f(E_{g_2}) = 4$.

|        | $E_{r_5}$ | $E_{r_1}$ | $E_{r_2}$ | $E_{r_3}$ | $E_{r_4}$ |
|--------|-----------|-----------|-----------|-----------|-----------|
|        |           |           |           | 1         | 4         |
|        |           |           | 2         |           |           |
|        |           |           |           | 4         | 1         |
| R, B1  | 2         | 1         |           |           |           |
|        |           |           |           | 1         | 4         |
|        |           |           | 4         |           |           |
|        |           |           |           | 2         | 1, 4      |

|        | $E_{r_5}$ | $E_{r_1}$ | $E_{r_2}$ | $E_{r_3}$ | $E_{r_4}$ |
|--------|-----------|-----------|-----------|-----------|-----------|
|        |           |           |           | 1         | 2         |
|        |           |           | 2         |           |           |
|        |           |           |           | 4         | 1, 2      |
| R, B2  | 4         | 1         |           |           |           |
|        |           |           |           | 1         | 2         |
|        |           |           | 4         |           |           |
|        |           |           |           | 2         | 1         |

in einer Tabelle

|  | $E_{r_1}$ | $E_{r_2}$ | $E_{r_3}$ | $E_{r_4}$ | $E_{r_5}$ |
|---|---|---|---|---|---|
|  | $E_{g_1}$ |  | $E_{g_4}$ | $E_{g_3}$ | $E_{g_2}$ |
|  |  |  | 4, 2 | 1 |  |
|  |  |  |  |  | 2 |
|  |  |  | 2 | 4 |  |
| G | 1 |  |  |  |  |
|  |  |  | 2, 4 | 1 |  |
|  |  |  |  |  | 4 |
|  |  |  | 4 | 2 |  |
| | | | 1 | 4 | |
|  |  | 2 |  |  |  |
|  |  |  | 4 | 1 |  |
| R, B1 | 1 |  |  |  | 2 |
|  |  |  | 1 | 4 |  |
|  |  | 4 |  |  |  |
|  |  |  | 2 | 1, 4 |  |
| | | | 1 | 2 | |
|  |  | 2 |  |  |  |
|  |  |  | 4 | 1, 2 |  |
| R, B2 | 1 |  |  |  | 4 |
|  |  |  | 1 | 2 |  |
|  |  | 4 |  |  |  |
|  |  |  | 2 | 1 |  |

finden wir als passende Färbungen

| $E_{r_1}$ | $E_{r_2}$ | $E_{r_3}$ | $E_{r_4}$ | $E_{r_5}$ |
|---|---|---|---|---|
| 1 | 2 | 4 | 1 | 2 |
|  | 4 | 2 | 1 | 2 |
|  | 4 | 2 | 4 | 2 |
|  | 4 | 2 | 1 | 4 |
|  | 2 | 4 | 1 | 4 |
|  | 2 | 4 | 2 | 4 |

als Block des FS(G+R)

| $E_{r_1}$ | $E_{r_2}$ | $E_{r_3}$ | $E_{r_4}$ |
|---|---|---|---|
|  |  | 2 | 4 |
|  | 1 |  |  |
|  |  | 4 | 2 |

Zu einem Wert $f(E_{g_1}) = 1$ gibt es die gekoppelten Spalten für $E_{r_{m-1}}$ $E_{r_m}$.

Dasselbe für $f(E_{g_1}) = 2$, $f(E_{g_2}) = 1, 4$.

| | $E_{r_1}$ | $E_{r_2}$ | $E_{r_3}$ | $E_{r_4}$ | $E_{r_5}$ |
|---|---|---|---|---|---|
| | $E_{g_1}$ | | $E_{g_4}$ | $E_{g_3}$ | $E_{g_2}$ |
| | | | 4 | 1 | |
| | | | | | 1 |
| | | | 2 | 4 | |
| G, B3 | 2 | | | | |
| | | | 4 | 2 | |
| | | | | | 4 |
| | | | 2, 4 | 1 | |
| | | | 2 | 4 | |
| | | 1 | | | |
| | | | 4 | 2 | |
| R, B3 | 2 | | | | 1 |
| | | | 2 | 4 | |
| | | 4 | | | |
| | | | 1 | 4 | |
| | | | 2 | 1 | |
| | | 1 | | | |
| | | | 4 | 2, 1 | |
| R, B4 | 2 | | | | 4 |
| | | | 2 | 1 | |
| | | 4 | | | |
| | | | 1 | 2 | |

finden wir als passende Färbungen

| | $E_{r_1}$ | $E_{r_2}$ | $E_{r_3}$ | $E_{r_4}$ | $E_{r_5}$ |
|---|---|---|---|---|---|
| | 2 | 1 | 4 | 2 | 1 |
| | | 4 | 2 | 4 | 1 |
| | | 1 | 2 | 4 | 1 |
| | | 1 | 4 | 2 | 4 |
| | | 1 | 2 | 1 | 4 |
| | | 4 | 2 | 1 | 4 |
| | | 1 | 4 | 1 | 4 |

als Block des FS(G+R)

| | $E_{r_1}$ | $E_{r_2}$ | $E_{r_3}$ | $E_{r_4}$ | $E_{r_5}$ |
|---|---|---|---|---|---|
| | | | 1 | 2, 4 | |
| | 2 | | | | |
| | | | 4 | 2 | |

| | $E_{r_1}$ | $E_{r_2}$ | $E_{r_3}$ | $E_{r_4}$ | $E_{r_5}$ |
|---|---|---|---|---|---|
| | $E_{g_1}$ | | $E_{g_4}$ | $E_{g_3}$ | $E_{g_2}$ |
| | | | 2 | 4 | |
| | | | | | 1 |
| | | | 4 | 2 | |
| G, B5 | 4 | | | | |
| | | | 2 | 4 | |
| | | | | | 2 |
| | | | 4, 2 | 1 | |

| | $E_{r_1}$ | $E_{r_2}$ | $E_{r_3}$ | $E_{r_4}$ | $E_{r_5}$ |
|---|---|---|---|---|---|
| | | | 4 | 2 | |
| | | 1 | | | |
| | | | 2 | 4 | |
| R, B5 | 4 | | | | 1 |
| | | | 4 | 2 | |
| | | 2 | | | |
| | | | 1 | 4, 2 | |

| | $E_{r_1}$ | $E_{r_2}$ | $E_{r_3}$ | $E_{r_4}$ | $E_{r_5}$ |
|---|---|---|---|---|---|
| | | | 4 | 1 | |
| | | 1 | | | |
| | | | 2 | 4, 1 | |
| R, B6 | 4 | | | | 2 |
| | | | 4 | 1 | |
| | | 2 | | | |
| | | | 1 | 4 | |

finden wir als passende Färbungen

| $E_{r_1}$ | $E_{r_2}$ | $E_{r_3}$ | $E_{r_4}$ | $E_{r_5}$ |
|---|---|---|---|---|
| 4 | 1 | 2 | 4 | 1 |
| 2 | 4 | 2 | | 1 |
| 1 | 2 | 4 | | 2 |
| 1 | 4 | 1 | | 2 |
| 2 | 4 | 1 | | 2 |
| 1 | 2 | 1 | | 2 |

als Block des FS(G+R)

| $E_{r_1}$ | $E_{r_2}$ | $E_{r_3}$ | $E_{r_4}$ | $E_{r_5}$ |
|---|---|---|---|---|
| | | 1 | 4, 2 | |
| | 4 | | | |
| | | 2 | 4 | |

2.2.3 k=4, m>4, $E_{g_1} = E_2$

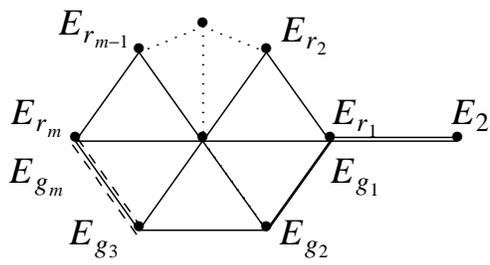

$f(E_{g_1}) = f(E_2) = 2.$

Wie oben, jedoch ist nur der Fall $f(E_{g_1}) = 2$ zu betrachten.

$E_{g_1} > E_2$

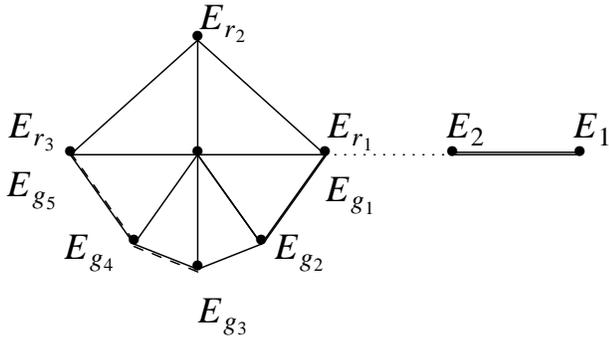

$E_{g_1} = E_2$

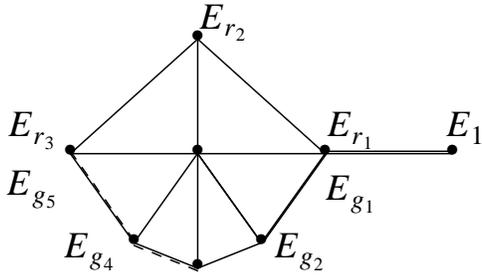

$f(E_{g_1}) = 1, 2$ oder 4, $f(E_{g_1}) = 2$ nur möglich, wenn $E_{g_1} = E_2$ oder wenn zwischen $E_{g_1} = E_2$ mindestens eine weitere Ecke liegt.

$f(E_{g_5}) \neq 1$ wegen $f(E_1) = 2$.

$f(E_{g_1}) = 1$, dann ist $f(E_{g_2}) = 2$ oder 4 und $f(E_{r_2}) = 2$ oder 4.

$f(E_{g_2}) = 2$.

$f(E_{g_5}) \neq 1 \Rightarrow$ G Block B3.

$f(E_{r_5}) \neq 2$ und $f(E_{r_2}) = 2 \Rightarrow$ R Block C3

$f(E_{r_5}) \neq 2$ und $f(E_{r_2}) = 4 \Rightarrow$ R Block B6

| | $E_{r_1}$ | $E_{r_2}$ | $E_{r_3}$ | $E_{r_4}$ | $E_{r_5}$ | $E_{g_2}$ |
| | $E_{g_1}$ | | $E_{g_5}$ | $E_{g_4}$ | $E_{g_3}$ | |
|---|---|---|---|---|---|---|
| | | | 4 | 2 | | |
| | | | | | 1 | |
| | | | 2 | 4 | | |
| G, B3 | 1 | | | | | 2 |
| | | | 4 | 2 | | |
| | | | | | 4 | |
| | | | 2, 4 | 1 | | |
| | | | | 2 | 4, 1 | |
| | | | 1 | | | |
| | | | | 4 | 1 | |
| R, C3 | 1 | 2 | | | | |
| | | | | 2 | 1, 4 | |
| | | | 4 | | | |
| | | | | 1 | 4 | |
| | | | | 4 | 1 | |
| | | | 1 | | | |
| | | | | 2 | 4, 1 | |
| R, B6 | 1 | 4 | | | | |
| | | | | 4 | 1 | |
| | | | 2 | | | |
| | | | | 1 | 4 | |

daraus die passenden Färbungen:

| | $E_{r_1}$ | $E_{r_2}$ | $E_{r_3}$ | $E_{r_4}$ | $E_{r_5}$ |
|---|---|---|---|---|---|
| | 1 | 2 | 4 | 2 | 1 |
| | | | 4 | 2 | 4 |
| | | | 4 | 1 | 4 |
| G+R | 1 | 4 | 2 | 4 | 1 |
| | | | 2 | 1 | 4 |

(benötigt wird natürlich nur eine passende Färbung um R anschließen zu können.)

also auf dem Rand von G+R

| $E_{r_1}$ | $E_{r_2}$ | $E_{r_3}$ |
|---|---|---|
| | 2 | 4 |
| 1 | | |
| | 4 | 2 |

$f(E_{g_2})$ ...

$f(E_{g_5}) \neq 1 \Rightarrow$ G Block Typ B5.

$f(E_{r_5}) \neq 4$ und $f(E_{r_2}) = 2 \Rightarrow$ R Block B4

$f(E_{r_5}) \neq 4$ und $f(E_{r_2}) = 4 \Rightarrow$ R Block C5

| | $E_{r_1}$ | $E_{r_2}$ | $E_{r_3}$ | $E_{r_4}$ | $E_{r_5}$ | |
| --- | --- | --- | --- | --- | --- | --- |
| | $E_{g_1}$ | | $E_{g_5}$ | $E_{g_4}$ | $E_{g_3}$ | $E_{g_2}$ |
| | | | 4 | 2 | | |
| | | | | | 1 | |
| | | | 2 | 4 | | |
| G, B5 | 1 | | | | | 4 |
| | | | 2 | 4 | | |
| | | | | 2 | | |
| | | | 2, 4 | 1 | | |

| | $E_{r_1}$ | $E_{r_2}$ | $E_{r_3}$ | $E_{r_4}$ | $E_{r_5}$ | |
| --- | --- | --- | --- | --- | --- | --- |
| | | | | 2 | 1 | |
| | | | 1 | | | |
| | | | | 4 | 2, 1 | |
| R, B4 | 1 | 2 | | | | |
| | | | | 2 | 1 | |
| | | | 4 | | | |
| | | | | 1 | 2 | |
| | | | | 4 | 2, 1 | |
| | | | 1 | | | |
| | | | | 2 | 1 | |
| R, C5 | 1 | 4 | | | | |
| | | | | 4 | 1, 2 | |
| | | | 2 | | | |
| | | | | 1 | 2 | |

daraus die passenden Färbungen:

| | $E_{r_1}$ | $E_{r_2}$ | $E_{r_3}$ | $E_{r_4}$ | $E_{r_5}$ | $E_{g_2}$ |
| --- | --- | --- | --- | --- | --- | --- |
| | | 1 | 2 | 4 | 2 | 1 |
| | | | | | 4 | 1 | 2 |
| G+R | | 1 | 4 | 2 | 4 | 1 |
| | | | | 2 | 4 | 2 |
| | | | | 2 | 1 | 2 |

| | $E_{r_2}$ | $E_{r_3}$ | $E_{r_4}$ | $E_{r_5}$ | $E_{g_2}$ |
| --- | --- | --- | --- | --- | --- |
| | 1 | 2 | 4 | 2 | 1 |
| | | | 4 | 1 | 2 |
| G+R | 1 | 4 | 2 | 4 | 1 |
| | | | 2 | 4 | 2 |
| | | | 2 | 1 | 2 |

also

| $E_{r_1}$ | $E_{r_2}$ | $E_{r_3}$ |
| --- | --- | --- |
| | 2 | 4 |
| 1 | | |
| | 4 | 2 |

$f(E_{g_1}) = 2$, dann ist $f(E_{g_2}) = 1$ oder 4 und $f(E_{r_2}) = 1$ oder 4.
$f(E_{g_2}) = 1$
$f(E_{g_5}) \neq 1 \Rightarrow$ G Block C1.
$f(E_{r_5}) \neq 1$ und $f(E_{r_2}) = 1 \Rightarrow$ R Block C1
$f(E_{r_5}) \neq 1$ und $f(E_{r_2}) = 4 \Rightarrow$ R Block B5

| | $E_{r_1}$ | $E_{r_2}$ | $E_{r_3}$ | $E_{r_4}$ | $E_{r_5}$ | |
|---|---|---|---|---|---|---|
| | $E_{g_1}$ | | $E_{g_5}$ | $E_{g_4}$ | $E_{g_3}$ | $E_{g_2}$ |
| | | | 4, 2 | 1 | | |
| | | | | | 2 | |
| | | | 2 | 4 | | |
| G, C1 | 2 | | | | | 1 |
| | | | 2, 4 | 1 | | |
| | | | | | 4 | |
| | | | 4 | 2 | | |
| | | | | 1 | 4, 2 | |
| | | | 2 | | | |
| | | | | 4 | 2 | |
| R, C1 | 2 | 1 | | | | |
| | | | | 1 | 2, 4 | |
| | | | 4 | | | |
| | | | | 2 | 4 | |
| | | | | 4 | 2 | |
| | | | 1 | | | |
| | | | | 2 | 4 | |
| R, B5 | 2 | 4 | | | | |
| | | | | 4 | 2 | |
| | | | 2 | | | |
| | | | | 1 | 4, 2 | |

daraus die passenden Färbungen:

| | | | | | |
|---|---|---|---|---|---|
| 2 | 1 | 4 | 1 | 2 |
| | | 2 | 1 | 2 |
| 2 | 4 | 2 | 1 | 2 |

also

| $E_{r_1}$ | $E_{r_2}$ | $E_{r_3}$ |
|---|---|---|
| | 1 | 4, 2 |
| 2 | | |
| | 4 | 2 |

$f(E_{g_1}) = 2, f(E_{g_2}) = 4$
$f(E_{g_5}) \neq 1 \Rightarrow$ G Block B5.
$f(E_{r_5}) \neq 4 \; f(E_{r_2}) = 1 \Rightarrow$ R Block B2
$f(E_{r_5}) \neq 4 \; f(E_{r_2}) = 4 \Rightarrow$ R Block C5

|        | $E_{r_1}$ | $E_{r_2}$ | $E_{r_3}$ | $E_{r_4}$ | $E_{r_5}$ | |
|--------|-----------|-----------|-----------|-----------|-----------|---|
|        | $E_{g_1}$ |           | $E_{g_5}$ | $E_{g_4}$ | $E_{g_3}$ | $E_{g_2}$ |
|        |           |           | 2         | 4         |           |   |
|        |           |           |           |           | 1         |   |
|        |           |           | 4         | 2         |           |   |
| G, B5  | 2         |           |           |           |           | 4 |
|        |           |           | 2         | 4         |           |   |
|        |           |           |           |           | 2         |   |
|        |           |           | 4, 2      | 1         |           |   |
| ─      |           |           |           |           |           |   |
|        |           |           |           | 1         | 2         |   |
|        |           |           | 2         |           |           |   |
|        |           |           |           | 4         | 1, 2      |   |
| R, B2  | 2         | 1         |           |           |           |   |
|        |           |           |           | 1         | 2         |   |
|        |           |           | 4         |           |           |   |
|        |           |           |           | 2         | 1         |   |
| ─      |           |           |           |           |           |   |
|        |           |           |           | 4         | 2, 1      |   |
|        |           |           | 1         |           |           |   |
|        |           |           |           | 2         | 1         |   |
| R, B5  | 2         | 4         |           |           |           |   |
|        |           |           |           | 4         | 1, 2      |   |
|        |           |           | 2         |           |           |   |
|        |           |           |           | 1         | 2         |   |

daraus die passenden Färbungen:

| $E_{r_1}$ | $E_{r_2}$ | $E_{r_3}$ | $E_{r_4}$ | $E_{r_5}$ |
|-----------|-----------|-----------|-----------|-----------|
| 2         | 1         | 2         | 4         | 1         |
|           |           | 4         | 2         | 1         |
|           | 2         | 4         | 2         |           |
| 2         | 4         | 2         | 4         | 1         |
|           |           | 2         | 4         | 2         |
|           |           | 2         | 1         | 2         |

also

| $E_{r_1}$ | $E_{r_2}$ | $E_{r_3}$ |
|-----------|-----------|-----------|
|           | 1         | 2, 4      |
| 2         |           |           |
|           | 4         | 2         |

$f(E_{g_2}) = 1$

$f(E_{g_5}) \neq 1 \Rightarrow$ G Block C1.

$f(E_{r_5}) \neq 1$ $f(E_{r_2}) = 1 \Rightarrow$ R Block C1.

$f(E_{r_5}) \neq 1$ $f(E_{r_2}) = 2 \Rightarrow$ R Block B3.

| | $E_{r_1}$ | $E_{r_2}$ | $E_{r_3}$ | $E_{r_4}$ | $E_{r_5}$ | |
| --- | --- | --- | --- | --- | --- | --- |
| | $E_{g_1}$ | | $E_{g_5}$ | $E_{g_4}$ | $E_{g_3}$ | $E_{g_2}$ |
| | | | 4, 2 | 1 | | |
| | | | | | 2 | |
| | | | 2 | 4 | | |
| G, C1 | 4 | | | | | 1 |
| | | | 4 | 2 | | |
| | | | | | 4 | |
| | | | 2, 4 | 1 | | |
| | | | 2 | 4, 1 | | |
| | | | 1 | | | |
| | | | 4 | 1 | | |
| R, B3 | 4 | 2 | | | | |
| | | | 1 | 4 | | |
| | | | 4 | | | |
| | | | 2 | 1, 4 | | |
| | | | 1 | 4, 2 | | |
| | | | 2 | | | |
| | | | 4 | 2 | | |
| R, C1 | | 1 | | | | |
| | | | 1 | 4 | | |
| | | | 4 | | | |
| | | | 2 | 4, 1 | | |

daraus die passenden Färbungen:

| $E_{r_1}$ | $E_{r_2}$ | $E_{r_3}$ | $E_{r_4}$ | $E_{r_5}$ |
| --- | --- | --- | --- | --- |
| 4 | 2 | 4 | 1 | 4 |
| | | 4 | 2 | 4 |
| | 1 | 2 | 1 | 2 |
| | | 2 | 4 | 2 |
| | | 2 | 1 | 4 |
| | | 4 | 1 | 4 |
| | | 4 | 2 | 1 |

also

| $E_{r_1}$ | $E_{r_2}$ | $E_{r_3}$ |
| --- | --- | --- |
| | 2 | 4 |
| 4 | | |
| | 1 | 2, 4 |

$f(E_{g_1}) = 4,$
$f(E_{g_2}) = 2 \Rightarrow$ G Typ B3.
$f(E_{r_5}) \neq 2, f(E_{r_2}) = 1 \Rightarrow$ R Typ B1.
$f(E_{r_5}) \neq 2, f(E_{r_2}) = 2 \Rightarrow$ G Typ C3.

| | $E_{r_1}$ | $E_{r_2}$ | $E_{r_3}$ | $E_{r_4}$ | $E_{r_5}$ | |
|---|---|---|---|---|---|---|
| | $E_{g_1}$ | | $E_{g_5}$ | $E_{g_4}$ | $E_{g_3}$ | $E_{g_2}$ |
| | | | 4 | 2 | | |
| | | | | | 1 | |
| | | | 2 | 4 | | |
| G, B3 | 4 | | | | | 2 |
| | | | 4 | 2 | | |
| | | | | | 4 | |
| | | | 2, 4 | 1 | | |

| | $E_{r_1}$ | $E_{r_2}$ | $E_{r_3}$ | $E_{r_4}$ | $E_{r_5}$ | |
|---|---|---|---|---|---|---|
| | | | | 1 | 4 | |
| | | | 2 | | | |
| | | | | 4 | 1 | |
| R, B1 | 4 | 1 | | | | |
| | | | | 1 | 4 | |
| | | | 4 | | | |
| | | | | 2 | 1, 4 | |

| | $E_{r_1}$ | $E_{r_2}$ | $E_{r_3}$ | $E_{r_4}$ | $E_{r_5}$ | |
|---|---|---|---|---|---|---|
| | | | | 2 | 4, 1 | |
| | | | 1 | | | |
| | | | | 4 | 1 | |
| R, C3 | 4 | 2 | | | | |
| | | | | 2 | 1, 4 | |
| | | | 4 | | | |
| | | | | 1 | 4 | |

daraus die passenden Färbungen:

| | $E_{r_1}$ | $E_{r_2}$ | $E_{r_3}$ | $E_{r_4}$ | $E_{r_5}$ |
|---|---|---|---|---|---|
| | 4 | 1 | 4 | 2 | 1 |
| | | | 2 | 4 | 1 |
| | | | 4 | 2 | 4 |
| | | | 2 | 1 | 4 |
| | | | 4 | 1 | 4 |
| | | 2 | 4 | 2 | 1 |
| | | | 4 | 2 | 4 |
| | | | 4 | 1 | 4 |

also

| $E_{r_1}$ | $E_{r_2}$ | $E_{r_3}$ |
|---|---|---|
| | 1 | 2, 4 |
| 1 | | |
| | 2 | 4 |

Die gekoppelten Ecken von G werden mit nicht gekoppelten Ecken von R identifiziert.

Die gekoppelten Ecken von R werden mit nicht gekoppelten Ecken von G identifiziert.

Die entsprechenden Blöcke des FS(G) zu den gekoppelten Ecken passen immer zu den normalen Blöcken des FS(R) der nicht gekoppelten Ecken.

Die entsprechenden Blöcke des FS(R) zu den gekoppelten Ecken passen immer zu den normalen Blöcken des FS(G) der nicht gekoppelten Ecken. (s. Tabelle oben.)

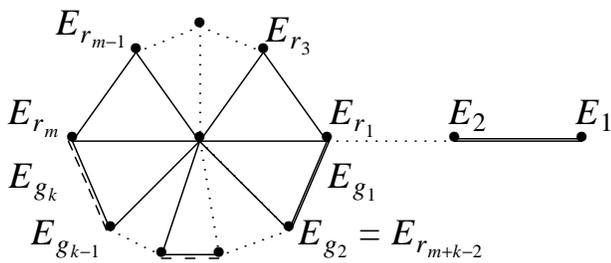

$E_{g_k}, E_{g_{k-1}}$ sind die gekoppelten Ecken von G.

$E_{r_{m+k-3}}, E_{r_{m+k-4}}$ sind die gekoppelten Ecken von R.

Die gekoppelten Ecken von G liegen nicht beide auf dem Rand von G+R.

Die gekoppelten Ecken von R liegen nicht auf dem Rand von G+R.

Die gekoppelten Ecken von G werden mit nicht gekoppelten Ecken von R identifiziert.

Die gekoppelten Ecken von R werden mit nicht gekoppelten Ecken von G identifiziert.

Bei der Konstruktion des FS(G+R) erwarten wir für $E_{r_{m+k-1}}, E_{g_{m+k}}$ gekoppelte Spalten.

Zwischen $E_{r_{k-1}}$ und $E_{g_4}$ werden Ecken identifiziert, deren Färbungen durch A1, A2 oder A3 beschrieben werden, also gibt es immer passende Färbungen: für $E_{g_{k-2}} \ E_{r_{m+2}}$, $E_{g_{k-3}} \ E_{r_{m+3}}$ , ..., gibt es immer passende Färbungen: zu $f(E_{g_{k-i}}) = f(E_{r_{m+i}}) = 1$ bzw. 2 bzw. 4 ist $f(E_{g_{k-i-1}}) = 2$ oder 4, bzw. 1 oder 4, bzw. 1 oder 2 und

$f(E_{r_{m+i+1}}) = 2$ und 4 bzw. 1 und 4 bzw. 1 und 2.

Wir betrachten

(1) die "Umgebung" von $E_{r_m}$

(2) die "Umgebung" von $E_{r_{m+k-4}}, E_{g_{m+k-3}} = E_{g_3}, E_{g_4}$.

Zu (1)

die möglichen Blöcke von G sind B3, B5, C1 ($f(E_{g_k}) \neq 1$),
die möglichen Blöcke von R sind A1, A2, A3, welcher hängt von m ab.
FS(G), B3 mit FS(R), A1

| | $E_{r_{m-2}}$ | $E_{r_{m-1}}$ | $E_{g_k}$ / $E_{r_m}$ | $E_{g_{k-1}}$ / $E_{r_{m+1}}$ | $E_{g_{k-2}}$ / $E_{r_{m+2}}$ | $E_{g_{k-3}}$ |
|---|---|---|---|---|---|---|
| | | | 4 | 2 | | |
| | | | | | 1 | |
| | | | 2 | 4 | | |
| G, B3 | | | 4 | 2 | | 2 |
| | | | | | 4 | |
| | | | 2, 4 | 1 | | |
| | | | | 2 | | |
| | | | 1 | | | |
| | | | | 4 | | |
| | | 2 | | 2 | | |
| | | | 4 | | | |
| | | | | 1 | | |
| R, A1 | 1 | | | 4 | | |
| | | | 1 | 2 | | |
| | | 4 | | 4 | | |
| | | | 2 | 1 | | |

die möglichen Färbungen für FS(G+R) als Block

| | $E_{r_{m-2}}$ | $E_{r_{m-1}}$ | $E_{g_k}$ / $E_{r_m}$ | $E_{g_{k-1}}$ / $E_{r_{m+1}}$ |
|---|---|---|---|---|
| | | | 2 | 4 |
| | 1 | | | |
| | | | 4 | 2 |

| | $E_{r_{m-2}}$ | $E_{r_{m-1}}$ | $E_{g_k}$ / $E_{r_m}$ | $E_{g_{k-1}}$ / $E_{r_{m+1}}$ | $E_{g_{k-2}}$ / $E_{r_{m+2}}$ | $E_{g_{k-3}}$ |
|---|---|---|---|---|---|---|
| | | | 4 | 2 | | |
| | | | | | 1 | |
| | | | 2 | 4 | | |
| G, B3 | | | | | | 2 |
| | | | 4 | 2 | | |
| | | | | | 4 | |
| | | | 2, 4 | 1 | | |

| | $E_{r_{m-2}}$ | $E_{r_{m-1}}$ | $E_{r_m}$ | $E_{r_{m+1}}$ | $E_{r_{m+2}}$ | |
|---|---|---|---|---|---|---|
| | | | | 1 | | |
| | | | 2 | | | |
| | | | | 4 | | |
| | | 1 | | 1 | | |
| | | | 4 | | | |
| | | | | 2 | | |
| R, A2 | 2 | | | | | |
| | | | | 4 | | |
| | | | 2 | | | |
| | | | | 1 | | |
| | | 4 | | | | |
| | | | | 4 | | |
| | | | 1 | | | |
| | | | | 2 | | |

die möglichen Färbungen für FS(G+R) als Block

| | $E_{r_{m-2}}$ | $E_{r_{m-1}}$ | $E_{r_m}$ | | | |
|---|---|---|---|---|---|---|
| | | 1 | 4, 2 | | | |
| | 2 | | | | | |
| | | 4 | 2 | | | |

| | $E_{r_{m-2}}$ | $E_{r_{m-1}}$ | $E_{r_m}$ | $E_{r_{m+1}}$ | | |
|---|---|---|---|---|---|---|
| | | | | 1 | | |
| | | | 4 | | | |
| | | | | 2 | | |
| | | 1 | | 1 | | |
| | | | 2 | | | |
| | | | | 4 | | |
| R, A3 | 4 | | | | | |
| | | | | 2 | | |
| | | | 4 | | | |
| | | | | 1 | | |
| | | 2 | | | | |
| | | | | 2 | | |
| | | | 1 | | | |
| | | | | 4 | | |

die möglichen Färbungen für FS(G+R) als Block

| | $E_{r_{m-2}}$ | $E_{r_{m-1}}$ | $E_{r_m}$ | | | |
|---|---|---|---|---|---|---|
| | | 1 | 2, 4 | | | |
| | 4 | | | | | |
| | | 2 | 4 | | | |

| | $E_{r_{m-2}}$ | $E_{r_{m-1}}$ | $E_{g_k}$ / $E_{r_m}$ | $E_{g_{k-1}}$ / $E_{r_{m+1}}$ | $E_{g_{k-2}}$ / $E_{r_{m+2}}$ | $E_{g_{k-3}}$ |
|---|---|---|---|---|---|---|
| | | | 2 | 4 | | |
| | | | | | 1 | |
| | | | 4 | 2 | | |
| G, B5 | | | 2 | 4 | | 4 |
| | | | | | 2 | |
| | | | 4, 2 | 1 | | |

| | $E_{r_{m-2}}$ | $E_{r_{m-1}}$ | $E_{g_k}$ / $E_{r_m}$ | $E_{g_{k-1}}$ / $E_{r_{m+1}}$ | $E_{g_{k-2}}$ / $E_{r_{m+2}}$ | $E_{g_{k-3}}$ |
|---|---|---|---|---|---|---|
| | | | | 2 | | |
| | | | 1 | | | |
| | | | | 4 | | |
| | | 2 | | | | |
| | | | | 2 | | |
| | | | 4 | | | |
| | | | | 1 | | |
| R, A1 | 1 | | | | | |
| | | | | 4 | | |
| | | | 1 | | | |
| | | | | 2 | | |
| | | 4 | | | | |
| | | | | 4 | | |
| | | | 2 | | | |
| | | | | 1 | | |
| | | | | 1 | | |
| | | | 2 | | | |

die möglichen Färbungen für FS(G+R) als Block

| | $E_{r_{m-2}}$ | $E_{r_{m-1}}$ | $E_{g_k}$ / $E_{r_m}$ | $E_{g_{k-1}}$ / $E_{r_{m+1}}$ | $E_{g_{k-2}}$ / $E_{r_{m+2}}$ | $E_{g_{k-3}}$ |
|---|---|---|---|---|---|---|
| | | | 2 | 4, 2 | | |
| | | 1 | | | | |
| | | | 2 | 4 | | |

| | $E_{r_{m-2}}$ | $E_{r_{m-1}}$ | $E_{g_k}$ / $E_{r_m}$ | $E_{g_{k-1}}$ / $E_{r_{m+1}}$ | $E_{g_{k-2}}$ / $E_{r_{m+2}}$ | $E_{g_{k-3}}$ |
|---|---|---|---|---|---|---|
| | | | | 4 | | |
| | | 1 | | | | |
| | | | | 1 | | |
| | | | 4 | | | |
| | | | | 2 | | |
| R, A2 | 2 | | | | | |
| | | | | 4 | | |
| | | | 2 | | | |
| | | | | 1 | | |
| | | 4 | | | | |
| | | | | 4 | | |
| | | | 1 | | | |
| | | | | 2 | | |

die möglichen Färbungen für FS(G+R) als Block

| | $E_{r_{m-2}}$ | $E_{r_{m-1}}$ | $E_{g_k}$ / $E_{r_m}$ | $E_{g_{k-1}}$ / $E_{r_{m+1}}$ | $E_{g_{k-2}}$ / $E_{r_{m+2}}$ | $E_{g_{k-3}}$ |
|---|---|---|---|---|---|---|
| | | 1 | 4, 2 | | | |
| | 2 | | | | | |
| | | 4 | 2 | | | |

| | $E_{r_{m-2}}$ | $E_{r_{m-1}}$ | $E_{g_k}$ / $E_{r_m}$ | $E_{g_{k-1}}$ / $E_{r_{m+1}}$ | $E_{g_{k-2}}$ / $E_{r_{m+2}}$ | $E_{g_{k-3}}$ |
|---|---|---|---|---|---|---|
| | | | 2 | 4 | | |
| | | | | | 1 | |
| | | | 4 | 2 | | |
| G, B5 | | | | | | 4 |
| | | | 2 | 4 | | |
| | | | | | 2 | |
| | | | 4, 2 | 1 | | |

| | $E_{r_{m-2}}$ | $E_{r_{m-1}}$ | $E_{r_m}$ | $E_{r_{m+1}}$ | $E_{r_{m+2}}$ | |
|---|---|---|---|---|---|---|
| | | | | 1 | | |
| | | | 4 | | | |
| | | | | 2 | | |
| | | 1 | | | | |
| | | | | 1 | | |
| | | | 2 | | | |
| | | | | 4 | | |
| R, A3 | 4 | | | | | |
| | | | | 2 | | |
| | | | 4 | | | |
| | | | | 1 | | |
| | | 2 | | | | |
| | | | | 2 | | |
| | | | 1 | | | |
| | | | | 4 | | |

die möglichen Färbungen für FS(G+R) als Block

| | | |
|---|---|---|
| | 1 | 2, 4 |
| 4 | | |
| | 2 | 4 |

FS(G), C1 mit FS(R), A1, A2, A3.

Wir erhalten völlig analog

| | | |
|---|---|---|
| | 2 | 4 |
| 1 | | |
| | 4 | 2 |
| | 1 | 4, 2 |
| 2 | | |
| | 4 | 2 |
| | 1 | 2, 4 |
| 4 | | |
| | 2 | 4 |

Zu (2)

Jetzt betrachten wir die gekoppelten Ecken von R, die mit normalen Ecken von G identifiert werden. Wir müssen lediglich zeigen, daß es passende Färbungen gibt.

Für $E_{g_1}$, $E_{g_2}$, $E_{g_3}$, $E_{g_4}$, haben wir im FS(G) A1, A2 oder A3, (falls $E_2 = E_{g_2}$ A2).
Mit $E_{g_1}$, $E_{g_2}$ als Basisseite von R und abhängig von $f(E_{g_2})$ sind für $E_{r_{m+k-6}}$, $E_{r_{m+k-5}}$, $E_{r_{m+k-4}}$, $E_{r_{m+k-3}}$ im FS(R) folgende Blöcke möglich:

| $f(E_{g_2})$ | Blöcke |
|---|---|
| 1 | B3 B5 C1 |
| 2 | B1 B6 C3 |
| 4 | B2 B4 C5 |

(Falls $E_2 = E_{g_1}$ ist $E_{r_1} = 2$ und es muß auch $f(E_{g_2}) = 1$ betrachtet werden.)

| Block G | $f(E_{g_2})$ | Blöcke von R |
|---|---|---|
| A1 | 2, 4 | B1 und B2 |
|  |  | B6 und C5 |
|  |  | C3 und B4 |
| A2 | 1, 4 | B3 und B4 |
|  |  | B5 und C5 |
|  |  | C1 und B2 |
| A2 | 1, 2 | B3 und C3 |
|  |  | B5 und B6 |
|  |  | C1 und B1 |

Wir betrachten alle Kombinationen für G, Block A1. Die Kombinationen für G, Block A2 und G, Block A3 gehen völlig analog.

$E_{r_{m+k-4}}$ .

| | $E_{r_{m+k-6}}$ | $E_{r_{m+k-5}}$ | $E_{g_4}$ / $E_{r_{m+k-4}}$ | $E_{g_3}$ / $E_{r_{m+k-3}}$ | $E_{g_2}$ | $E_{g_1}$ |
|---|---|---|---|---|---|---|
| | | | 2 | 1 | | |
| | | | 4 | | 2 | |
| | | | 2 | 4 | | |
| G, A1 | | | 1 | | | 1 |
| | | | 4 | 1 | | |
| | | | 2 | | 4 | |
| | | | 4 | 2 | | |
| | | | 1 | | | |
| | | | 1 | 4 | | |
| | | 2 | | | | |
| | | | 4 | 1 | | |
| R, B1 | 1 | | | | 2 | |
| | | | 1 | 4 | | |
| | | 4 | | | | |
| | | | 2 | 1, 4 | | |
| | | | 1 | 2 | | |
| | | 2 | | | | |
| | | | 4 | 1, 2 | | |
| R, B2 | 1 | | | | 4 | |
| | | | 1 | 2 | | |
| | | 4 | | | | |
| | | | 2 | 1 | | |

Wir stellen fest, daß es zu jeder Färbung im FS(G) eine passende Färbung in B1 oder B2 gibt und umgekehrt.

| | $E_{r_{m+k-6}}$ | $E_{r_{m+k-5}}$ | $E_{r_{m+k-4}}$ | $E_{r_{m+k-3}}$ | |
|---|---|---|---|---|---|
| | | | 2 | | |
| | | | | 1 | |
| | | | 4 | | |
| | | | | | 2 |
| | | | 2 | | |
| | | | | 4 | |
| G, A1 | | | 1 | | 1 |
| | | | 4 | | |
| | | | | 1 | |
| | | | 2 | | |
| | | | | | 4 |
| | | | 4 | | |
| | | | | 2 | |
| | | | 1 | | |
| | | | 4 | 1 | |
| | | 2 | | | |
| | | | 1 | 4 | |
| R, B6 | 4 | | | | 2 |
| | | | 4 | 1 | |
| | | 1 | | | |
| | | | 2 | 4, 1 | |
| | | | 4 | 2, 1 | |
| | | 1 | | | |
| | | | 2 | 1 | |
| R, C3 | 4 | | | | 4 |
| | | | 4 | 1, 2 | |
| | | 2 | | | |
| | | | 1 | 2 | |

Wir stellen fest, daß es zu jeder Färbung im FS(G) eine passende Färbung in B6 oder C3 gibt und umgekehrt.

| | $E_{r_{m+k-6}}$ | $E_{r_{m+k-5}}$ | $E_{r_{m+k-4}}$ | $E_{r_{m+k-3}}$ | |
|---|---|---|---|---|---|
| | | | 2 | | |
| | | | | 1 | |
| | | | 4 | | |
| | | | | | 2 |
| | | | 2 | 4 | |
| | | | 1 | | |
| G, A1 | | | | | 1 |
| | | | 4 | 1 | |
| | | | 2 | | 4 |
| | | | 4 | 2 | |
| | | | 1 | | |
| | | | 2 | 4, 1 | |
| | | 1 | | | |
| | | | 4 | 1 | |
| R, C2 | 2 | | | | 2 |
| | | | 2 | 1, 4 | |
| | | 4 | | | |
| | | | 1 | 4 | |
| | | | 2 | 1 | |
| | | 1 | | | |
| | | | 4 | 2, 1 | |
| R, B4 | 2 | | | | 4 |
| | | | 1 | 2 | |
| | | 4 | | | |
| | | | 2 | 1 | |

Wir stellen fest, daß es zu jeder Färbung im FS(G) eine passende Färbung in C2 oder B4 gibt und umgekehrt.

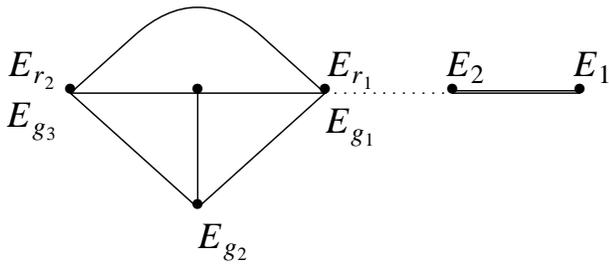

$E_{r_2}$ $E_{r_1}$ $E_2$ $E_1$
$E_{g_3}$ $E_{g_1}$
$E_{g_2}$

Wir betrachten die möglichen Färbungen von $E_{g_1}$ und die entsprechenden Blöcke des FS(G). Daraus ergeben sich die Färbungen von $E_{r_2}$.

(R ist maximal, $\gamma(R) = 3$, hat also passende Färbungen, wenn $f(E_{g_1}) \neq f(E_{g_2}) \neq f(E_{g_3}) \neq f(E_{g_1})$.)

|  | $E_{g_1}$ | $E_{g_2}$ | $E_{g_3}$ |  |
|---|---|---|---|---|
|  |  | 2 | 4 |  |
| $f(E_{g_1}) = 1 => f(E_{g_2}) = 2,4$   G | 1 |  |  | $f(E_{r_2}) = 2,4$ |
|  |  | 4 | 2 |  |
|  |  | 1 | 4,2 |  |
| $f(E_{g_1}) = 2 => f(E_{g_2}) = 1,4$   G | 2 |  |  | $f(E_{r_2}) = 4$ |
|  |  | 4 | 2 |  |
|  |  | 1 | 2,4 |  |
| $f(E_{g_1}) = 4 => f(E_{g_2}) = 1,2$   G | 4 |  |  | $f(E_{r_2}) = 2$ |
|  |  | 4 | 2 |  |

Nun betrachten wir $E_{g_1-1}$ (falls $E_{g_1-1} = E_2$ ist nur $f(E_{g_1-1}) = 2$ zutreffend) und erhalten die Blöcke des FS(G+R):

|  | $E_{g_1-1}$ | $E_{r_1}$ | $E_{r_2}$ |
|---|---|---|---|
|  |  | 2 | 4 |
| $f(E_{g_1-1}) = 1 => f(E_{g_1}) = 2,4$   G+R | 1 |  |  |
|  |  | 4 | 2 |
|  |  | 1 | 4,2 |
| $f(E_{g_1-1}) = 2 => f(E_{g_1}) = 1,4$   G+R | 2 |  |  |
|  |  | 4 | 2 |
|  |  | 1 | 2,4 |
| $f(E_{g_1-1}) = 4 => f(E_{g_1}) = 1,2$   G+R | 4 |  |  |
|  |  | 2 | 4 |

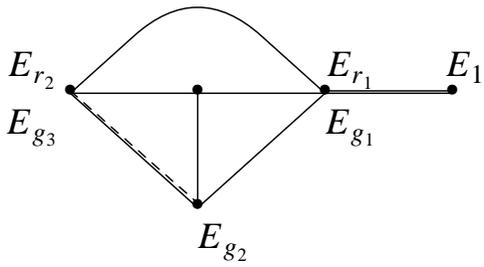

$\gamma(G) = 4.$

|        | $E_1$ | $E_{g_1}$ | $E_{g_2}$ | $E_{g_3}$ |
|--------|-------|-----------|-----------|-----------|
|        |       |           | 1         | 4, 2      |
| FS(G) = | 1    | 2         |           |           |
|        |       |           | 4         | 2         |

also ist

|          | $E_1$ | $E_{r_1}$ | $E_{r_2}$ |
|----------|-------|-----------|-----------|
| FS(G+R) = | 1    | 2         | 4         |

(3.2.1 nächste Seite)
3.2.2 k=4, m=2, $E_{g_1} = E_2$

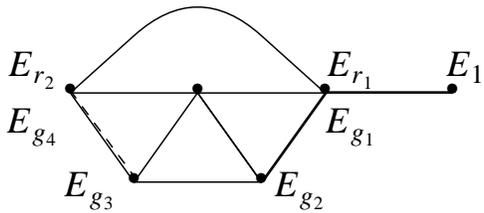

$\gamma(G) = 5.$

|        | $E_1$ | $E_{g_1}$ | $E_{g_2}$ | $E_{g_3}$ | $E_{g_4}$ |
|--------|-------|-----------|-----------|-----------|-----------|
|        |       |           |           | 2         | 4         |
|        |       |           | 1         |           |           |
|        |       |           |           | 4         | 2         |
| FS(G) = | 1    | 2         |           |           |           |
|        |       |           |           | 2         | 4         |
|        |       |           | 4         |           |           |
|        |       |           |           | 1         | 2, 4      |

also ist

|          | $E_1$ | $E_{r_1}$ | $E_{r_2}$ |
|----------|-------|-----------|-----------|
| FS(G+R) = | 1    | 2         | 4         |

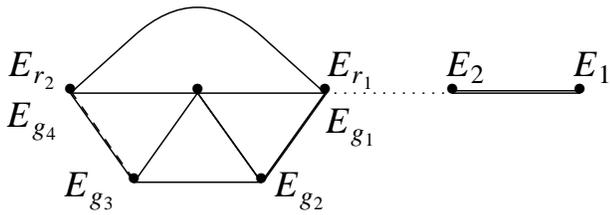

Wir betrachten die möglichen Färbungen von $E_{g_1}$ und die entsprechenden Blöcke des FS(G). Daraus ergeben sich die Färbungen von $E_{r_2}$.

Sei $f(E_{g_1}) = 1$. (Block von FS(G) ist C1.)

| $E_{r_1}$ | $E_{r_2}$ | $E_{r_3}$ | $E_{r_4}$ |
|-----------|-----------|-----------|-----------|
| $E_{g_1}$ | $E_{g_4}$ | $E_{g_3}$ | $E_{g_2}$ |
|           | 4, 2      | 1         |           |
|           |           |           | 2         |
|           | 2         | 4         |           |
| 1         |           |           |           |
|           | 2, 4      | 1         |           |
|           |           |           | 4         |
|           | 4         | 2         |           |

FS(R), $f(E_{r_4}) = 2$

|   | 2 | 1, 4 |   |
|---|---|------|---|
|   |   |      | 2 |
|   | 4 | 1    |   |

FS(R), $f(E_{r_4}) = 4$

|   | 2 | 1    |   |
|---|---|------|---|
| 1 |   |      | 4 |
|   | 4 | 2, 1 |   |

passende Färbungen

| 1 | 4 | 1 | 2 |
|---|---|---|---|
|   | 2 | 1 | 2 |
|   | 2 | 4 | 2 |
|   | 2 | 1 | 4 |
|   | 4 | 1 | 4 |
|   | 4 | 2 | 4 |

auf dem Rand von G+R

| 1 | 2, 4 |
|---|------|

Sei $f(E_{g_1}) = 2$.

Dies ist der einzige zu betrachtende Fall, wenn $E_{g_1} = E_2$.

| $E_{r_1}$ | $E_{r_2}$ | $E_{r_3}$ | $E_{r_4}$ |
|---|---|---|---|
| $E_{g_1}$ | $E_{g_4}$ | $E_{g_3}$ | $E_{g_2}$ |
|  | 4 | 2 |  |
|  |  |  | 1 |
|  | 2 | 4 |  |
| 2 |  |  |  |
|  | 4 | 2 |  |
|  |  |  | 4 |
|  | 2, 4 | 1 |  |

FS(R), $f(E_{r_4}) = 1$

| | | | |
|---|---|---|---|
|  | 1 | 4, 2 |  |
| 2 |  |  | 1 |
|  | 4 | 2 |  |

FS(R), $f(E_{r_4}) = 4$

| | | | |
|---|---|---|---|
|  | 1 | 2 |  |
| 2 |  |  | 4 |
|  | 4 | 1, 2 |  |

passende Färbungen

| | | | |
|---|---|---|---|
| 2 | 4 | 2 | 1 |
|  | 4 | 2 | 4 |
|  | 4 | 1 | 4 |

auf dem Rand von G+R

| | | | |
|---|---|---|---|
| 2 | 4 |  |  |

Sei $f(E_{g_1}) = 4$.

| $E_{r_1}$ | $E_{r_2}$ | $E_{r_3}$ | $E_{r_4}$ |
|---|---|---|---|
| $E_{g_1}$ | $E_{g_4}$ | $E_{g_3}$ | $E_{g_2}$ |
|  | 2 | 4 |  |
|  |  |  | 1 |
|  | 4 | 2 |  |
| 4 |  |  |  |
|  | 2 | 4 |  |
|  |  |  | 2 |
|  | 2, 4 | 1 |  |

FS(R), $f(E_{r_4}) = 1$

| $E_{r_1}$ | $E_{r_2}$ | $E_{r_3}$ | $E_{r_4}$ |
|---|---|---|---|
|  | 1 | 4, 2 |  |
| 4 |  |  | 1 |
|  | 2 | 4 |  |

FS(R), $f(E_{r_4}) = 2$

| $E_{r_1}$ | $E_{r_2}$ | $E_{r_3}$ | $E_{r_4}$ |
|---|---|---|---|
|  | 1 | 4 |  |
| 4 |  |  | 2 |
|  | 2 | 1, 4 |  |

passende Färbungen

| $E_{r_1}$ | $E_{r_2}$ | $E_{r_3}$ | $E_{r_4}$ |
|---|---|---|---|
| 4 | 2 | 4 | 1 |
|  | 2 | 4 | 2 |
|  | 2 | 1 | 2 |

auf dem Rand von G+R

| $E_{r_1}$ | $E_{r_2}$ |
|---|---|
| 4 | 2 |

Für FS(G+R) erhalten wir

| $f(E_{g_1-1})$ | $f(E_{g_1})$ | $E_{g_1-1}$ | $E_{r_1}$ | $E_{r_2}$ |
|---|---|---|---|---|
|  |  |  | 2 | 4 |
| 1 | 2, 4 | 1 |  |  |
|  |  |  | 4 | 2 |
|  |  |  | 1 | 4, 2 |
| 2 | 1, 4 | 2 |  |  |
|  |  |  | 4 | 2 |
|  |  |  | 1 | 2, 4 |
| 4 | 1, 2 | 4 |  |  |
|  |  |  | 2 | 4 |

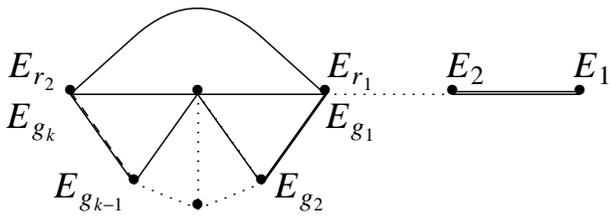

R anschließen mit der gleichen Orientierung.

Zu jedem Element der Spalte von $E_{g_1}$ gehört ein TeilFS von G das wir als FS(R) betrachten können. $E_{g_1}$ $E_{g_2}$ als Basisseite. Wegen der Kante $E_{r_1}$ $E_{r_2}$ fallen in diesem FS in der Spalte von $E_{r_2} = E_{g_k}$ alle Elemente weg, die gleich $f(E_{g_1})$ sind.

Wir betrachten FS(G), der letzte Block kann sein: B3, B5, C1

Wir betrachten B3 und die möglichen Färbungen von $E_{r_2}$ in Abhängigkeit von $f(E_{g_1})$. ($E_{r_i} = E_{g_{k+2-i}}$, i=2, ..., 5)

| $E_{g_k}$ | $E_{g_{k-1}}$ | $E_{g_{k-2}}$ | $E_{g_{k-3}}$ | $E_{g_1}$ | $E_{r_2}$ |
|---|---|---|---|---|---|
|  |  | 2 | 4 | 1 | 2, 4 |
|  | 1 |  |  | 2 | 4 |
|  |  | 4 | 2 | 4 | 2 |
| 2 |  |  |  |  |  |
|  |  | 2 | 4 |  |  |
|  | 4 |  |  |  |  |
|  |  | 1 | 2, 4 |  |  |

Dasselbe für B5

| $E_{g_k}$ | $E_{g_{k-1}}$ | $E_{g_{k-2}}$ | $E_{g_{k-3}}$ | $E_{g_1}$ | $E_{r_2}$ |
|---|---|---|---|---|---|
|  |  | 4 | 2 | 1 | 2, 4 |
|  | 1 |  |  | 2 | 4 |
|  |  | 2 | 4 | 4 | 2 |
| 4 |  |  |  |  |  |
|  |  | 4 | 2 |  |  |
|  | 2 |  |  |  |  |
|  |  | 1 | 4, 2 |  |  |

und für C1

| $E_{g_k}$ | $E_{g_{k-1}}$ | $E_{g_{k-2}}$ | $E_{g_{k-3}}$ | $E_{g_1}$ | $E_{r_2}$ |
|---|---|---|---|---|---|
|  |  | 1 | 4, 2 | 1 | 4, 2 |
|  | 2 |  |  | 2 | 4 |
|  |  | 4 | 2 | 4 | 2 |
| 1 |  |  |  |  |  |
|  |  | 1 | 2, 4 |  |  |
|  | 4 |  |  |  |  |
|  |  | 2 | 4 |  |  |

Betrachtung dazu die Werte für $f(E_{r_2})$, dargestellt als Blöcke:

| $E_{g_1-1}$ | $E_{g_1}$ | $E_{r_2}$ |
|---|---|---|
| | 2 | 4 |
| 1 | | |
| | 4 | 2 |
| | | |
| | 1 | 2, 4 |
| 2 | | |
| | 4 | 2 |
| | | |
| | 1 | 4, 2 |
| 4 | | |
| | 2 | 4 |

Wir erhalten FS(G+R) mit $E_{r_1}$ $E_{r_2}$ als gekoppelte Spalten. Falls $E_{g_1-1} = E_2$ ist nur $f(E_{g_1-1}) = 2$ möglich und es ist $\gamma(G + R) = 4$, dann erhalten wir:

| | | | 1 | 2, 4 |
|---|---|---|---|---|
| FS(G+R) = | 1 | 2 | | |
| | | | 4 | 2 |

Wenn $E_{g_1} = E_2$ dann ist $f(E_{g_1}) = 2$ und $f(E_{r_1}) = 4$ und G+R ist maximal, FS(G+R) = 1 2 4.

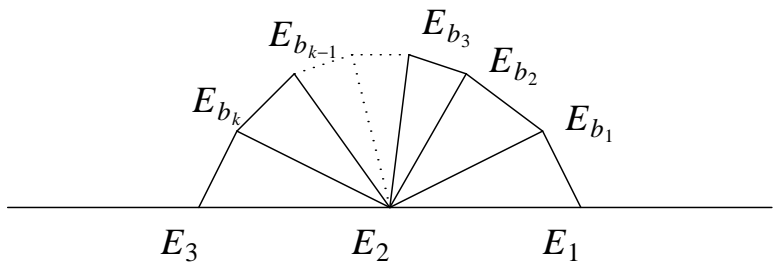

Fig Bogen

Dargestellt ist ein Bogen der Länge k.

Gegeben sei ein fm Graph G mit Färbungsschema FS(G). Der Bogen wird wie dargestellt an G angeschlossen.

Wir betrachten einen Block zu einem Element der Spalte von $E_1$:

| $E_1$ | $E_2$ | $E_3$ |
|---|---|---|
|  |  | 1 |
|  | 2 |  |
|  |  | 4 |
|  |  |  |
| 1 |  |  |
|  |  |  |
|  |  | 1 |
|  | 4 |  |
|  |  | 2 |

Wir zeigen nun die FS zu den Graphen, die wir erhalten indem wir die Ecken des Bogens der Reihe nach anschliessen, also $G + E_{b_1}$, $G + E_{b_1} + E_{b_2}$, ...

| $E_1$ | $E_{b_1}$ | $E_2$ | $E_3$ | als Liste | | | |
|---|---|---|---|---|---|---|---|
|  | 3 |  | 1 | 1 | 3 | 2 | 1 |
|  |  | 2 |  | 1 | 4 | 2 | 1 |
|  | 4 |  | 4 | 1 | 3 | 2 | 4 |
|  |  |  |  | 1 | 4 | 2 | 4 |
| 1 |  |  |  | 1 | 3 | 4 | 1 |
|  | 3 |  | 1 | 1 | 2 | 4 | 1 |
|  |  | 4 |  | 1 | 3 | 4 | 2 |
|  | 2 |  | 2 | 1 | 2 | 4 | 2 |

Nach dem Anschluß der zweiten Ecke des Bogens haben wir als FS

| $E_1$ | $E_{b_1}$ | $E_{b_2}$ | $E_2$ | $E_3$ |
|---|---|---|---|---|
|  |  | 1 |  |  |
|  | 3 |  |  | 1 |
|  |  | 4 |  |  |
|  |  |  | 2 |  |
|  |  | 1 |  |  |
|  | 4 |  |  | 4 |
|  |  | 3 |  |  |
| 1 |  |  |  |  |
|  |  | ... |  |  |
|  |  |  | ... |  |
|  |  | ... |  |  |

Die Darstellung der unteren Hälfte können wir uns sparen (vertausche 2 und 4).

```
       1
   3       1
       1
   4
       3
         2
       1
   3
       4
       1
   4       4
       3
1
```

Zwischen den Spalten von $E_1$ und $E_2$ sehen wir jetzt ein "ganz normales" FS, sozusagen ein lokales FS, wobei die Ecken $E_1$ und $E_2$ mit ihren Farben die analoge Rolle spielen wie die Farben der Ecken $E_2$ und $E_3$ der Basisseite.

Interessant wird der Anschluß der letzten Ecke des Bogens, da für diese die zulässigen Farben von den Farben von drei Ecken abhängen. Wir erhalten - wie zu erwarten - gekoppelte Spalten.

Die behauptete Anzahl auf dem Rand verschiedener Färbungen ist gegeben.
Es handelt sich gemäß Konstruktion um lauter zulässige Färbungen.

In diesem FS kommen als Blöcke nicht nur die bekannten Typen vor. Solche FS würden den weiteren Beweis komplizierter machen, deshalb erfolgt der Anschluß eines Bogens nur wenn die Ecken des Bogens auf dem Rand des letztendlich zu färbenden Graphen G liegen, d.h. wenn keine Ecke des Bogens durch den Anschluß weiterer Ecken nach innen fällt.

Bemerkung:

Natürlich ist ein Bogen ein Rad mit $E_2$ als Nabe. Die Färbungen für $E_2$ stehen aber schon fest. Um statt des Bogens ein Rad mit $E_2$ als Nabe anschließen zu können, müßte man die Konstruktion des Graphen zum Teil rückgängig machen (soweit, daß $E_2$ noch nicht angeschlossen ist), dies möchten wir vermeiden.

Bei mehreren Bögen auf dem Rand gibt es nichts prinzipiell Neues wenn sich die Bögen nicht überlappen.

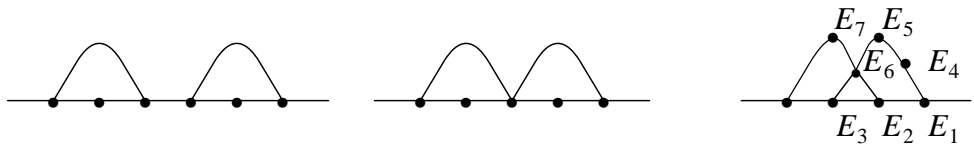

Im dritten Beispiel überlappen sich zwei Bögen. Auch hier können wir die behaupteten Färbungen konstruieren:

| $E_1$ | $E_4$ | $E_5$ | $E_2$ | $E_3$ | $E_1$ | $E_4$ | $E_5$ | $E_6$ | $E_7$ | $E_3$ |
|---|---|---|---|---|---|---|---|---|---|---|
| | | 1 | | | | | 1 | 4, 3 | | |
| | 3 | | | 1 | | 3 | | | | 1 |
| | | 4 | | | | | 4 | 3 | | |
| | | | 2 | | | | | | | |
| | | 1 | | | | | 1 | 3 | | |
| | 4 | | | 4 | | 4 | | | | 4 |
| | | 3 | | | | | 3 | 1 | | |
| 1 | | | | | 1 | | | | | |
| | | 2 | | | | | 2 | 3 | | |
| | 3 | | | 1 | | 3 | | | | 1 |
| | | 1 | | | | | 1 | 2, 3 | | |
| | | | 4 | | | | | | | |
| | | 2 | | | | | 2 | 3, 1 | | |
| | 1 | | | 2 | | 1 | | | | 2 |
| | | 3 | | | | | 3 | 1 | | |

Die Spalte für $E_7$ überlassen wir dem Leser. (Danke)

Das FS wird wieder kompliziert und diese Konstruktion wird nur verwendet, wenn die Ecken der Bögen auf dem Rand (von G) bleiben. (Wenn sie durch weitere Anschlüsse nach innen fallen ist der Anschluß eines Rades möglich.)

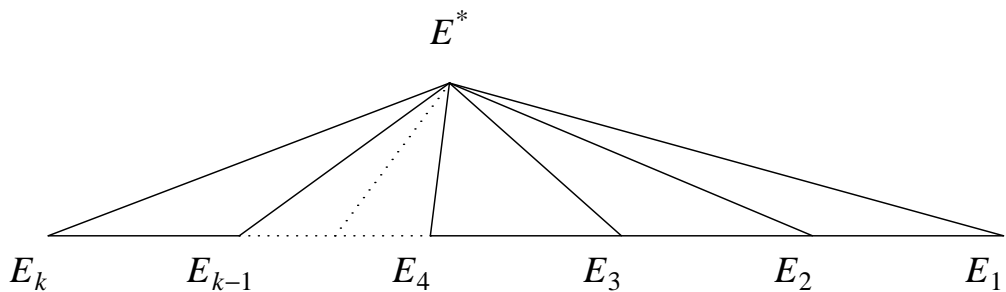

$$E^*$$

$$E_k \quad E_{k-1} \quad E_4 \quad E_3 \quad E_2 \quad E_1$$

Dargestellt ist eine Ecke $E^*$, die durch Kanten mit k Ecken eines fm Graphen verbunden wird.

Es interessiert nur $k \geq 3$; der Fall $k = 2$ ist der Anfang eines Bogens.

Dieser Anschluß ist nur zu betrachten für den Fall, daß $E^*$ auf dem Rand des zu färbenden Graphen liegt, d.h. nicht nach weiteren Anschlüssen nach innen fällt.

Wir betrachten das FS(G), und zwar die drei Spalten zu den Ecken $E_1, E_2, E_3$, davon einen Block. Für alle anderen Blöcke verläuft die Konstruktion der Färbungen analog.

Im allgemeinen sind die Spalten der Ecken $E_1$ usw. länger, aber alle Blöcke zu jeweils 4 Färbungen sehen im Prinzip gleich aus.

| $E_1$ | $E_2$ | $E_3$ | | $E_1$ | $E^*$ | $E_3$ |
|-------|-------|-------|--|-------|-------|-------|
|       |       | 1     |  |       | 4, 3  | 1     |
|       | 2     |       |  |       |       |       |
|       |       | 4     |  |       | 3     | 4     |
|       |       |       |  |       |       |       |
| 1     |       |       | wird zu | 1 | | |
|       |       |       |  |       |       |       |
|       |       | 1     |  |       | 2, 3  | 1     |
|       | 4     |       |  |       |       |       |
|       |       | 2     |  |       | 3     | 2     |

$E^*$ hat also die gleiche Länge wie $E_3$. In der Spalte jeder weiteren Ecke $(E_4, \ldots, E_k)$, die mit $E^*$ verbunden wird, werden die Farben, die wegen der jeweiligen Farbe von $E^*$ nicht zulässig (da gleich) sind, gestrichen, also maximal jede zweite. Dadurch hat diese Ecke dann die Länge von $E^*$ und dabei fallen $E_4, \ldots, E_{k-1}$ nach innen. Wir haben schließlich die behaupteten Färbungen (mindestens) für den Rand von $G + E^*$ mit zwei gekoppelten (d. h. gleich langen) Spalten für $E^*, E_k$.

Bemerkung: Obiges FS in der Mitte legt es nahe für $f(E^*) = 3$ zu wählen und zu zeigen, daß es die behauptete Anzahl auf dem Rand von $G + E^*$ sich unterscheidende Färbungen gibt. Dies ist aber nicht unbedingt richtig: Zeile 1 und Zeile 3 sind mit dieser Färbung für $E^*$ (zunächst lokal) gleich und daher kann es Färbungen geben, die auf dem ganzen Rand von $G + E^*$ gleich sind.

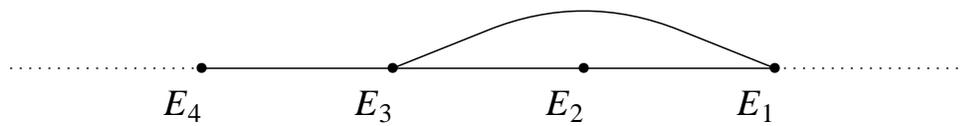

$$E_4 \qquad E_3 \qquad E_2 \qquad E_1$$

Fig Kante

Gegeben sei ein fm Graph G mit Färbungsschema FS. Eine Kante K werde angefügt. Damit G' = G + K fm ist, muß K eine Ecke des Randes von G mit der übernächsten verbinden. Eine Kante von $E_1$ nach $E_4$ ohne die eingezeichnete Kante gibt einen Graph, der nicht fm ist.

Wir betrachten den Fall, daß beim Anschluß von K keine gekoppelte Spalte nach innen fällt. Das FS von G

| $E_1$ | $E_2$ | $E_3$ | |
|---|---|---|---|
| | | 1 | < |
| | 2 | | |
| | | 4 | < |
| | | | |
| 1 | | | |
| | | | |
| | | 1 | < |
| | 4 | | |
| | | 2 | < |

wird zum FS von G'

| $E_1$ | $E_3$ | |
|---|---|---|
| | 4 | < |
| | | |
| 1 | | |
| | | |
| | 2 | < |

Dabei deuten die < die Fortsetzung des FS an, d.h. die weiteren Spalten, an denen sich nichts ändert.

Die Spalte von $E_2$, Spalte zwei im FS(G), taucht im FS(G') nicht mehr auf, da die Ecke nach innen fällt. In der dritten Spalte von FS(G) fallen die Färbungen mit "1" weg, da sie jetzt unzulässig sind, d.h. die Anzahl der Färbungen halbiert sich.

Nun betrachten wir den Fall, daß beim Anschluß von K eine gekoppelte Spalte nach innen fällt. Dazu betrachten wir die FS(R5) und FS(R6), Abschnitt 5. Färbungen Räder. (Eckennumerierung nicht die der Fig Kante.)

| $E_1$ | $E_2$ | $E_4$ | $E_5$ | $E_6$ |
|-------|-------|-------|-------|-------|
|       |       |       | 2     | 4     |
|       |       | 1     |       |       |
|       |       |       | 4     | 2     |
|       |       |       |       |       |
| 1     | 2     |       |       |       |
|       |       |       |       |       |
|       |       |       | 2     | 4     |
|       |       | 4     |       |       |
|       |       |       | 1     | 2, 4  |

Die Kante verbinde nun $E_4$ mit $E_6$. Wir erhalten als FS:

| $E_1$ | $E_2$ | $E_4$ | $E_6$ |
|-------|-------|-------|-------|
|       |       | 1     | 4, 2  |
|       |       |       |       |
| 1     | 2     |       |       |
|       |       |       |       |
|       |       | 4     | 2     |

Entsprechend für R6 mit der Kante von $E_5$ nach $E_7$:

| $E_1$ | $E_2$ | $E_4$ | $E_5$ | $E_6$ | $E_7$ |
|-------|-------|-------|-------|-------|-------|
|       |       |       |       | 1     | 4, 2  |
|       |       |       | 2     |       |       |
|       |       |       |       | 4     | 2     |
|       |       | 1     |       |       |       |
|       |       |       |       | 1     | 2, 4  |
|       |       |       | 4     |       |       |
|       |       |       |       | 2     | 4     |
|       |       |       |       |       |       |
| 1     | 2     |       |       |       |       |
|       |       |       |       | 4     | 2     |
|       |       |       | 2     |       |       |
|       |       |       |       | 1     | 4, 2  |
|       |       | 4     |       |       |       |
|       |       |       |       | 4     | 2     |
|       |       |       | 1     |       |       |
|       |       |       |       | 2     | 4     |

Wir erhalten das FS:

| $E_1$ | $E_2$ | $E_4$ | $E_5$ | $E_7$ |
|-------|-------|-------|-------|-------|
|       |       |       | 2     | 4     |
|       |       | 1     |       |       |
|       |       |       | 4     | 2     |
|       |       |       |       |       |
| 1     | 2     |       |       |       |
|       |       |       | 2     | 4     |
|       |       | 4     |       |       |
|       |       |       | 1     | 2, 4  |

Wir erhalten also wieder ein FS mit gekoppelten Spalten. Der Länge des Randes ist um eins kleiner und die Anzahl der Färbungen ist halbiert.

Tatsächlich haben uns nur die Blöcke interessiert. Das Ergebnis ist, daß wir - nach dem Anschluß der Kante - wieder ein FS haben in dem nur die bekannten Blöcke vorkommen.

## 11. Spezialfälle

Einfache Beispiele zeigen, daß nicht alle plättbare Graphen zu K gehören. Wir betrachten Beispiele.

Beispiel 1

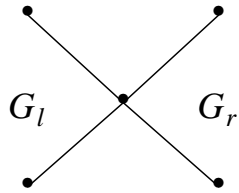

Es seien $G_l \in K$ und $G_r \in K$, mit einer gemeinsamen Ecke bilden sie G (G ist kein i-Graph). Mit $g = \gamma(G)$, $g_l = \gamma(G_l)$, $g_r = \gamma(G_r)$ ist $g = g_l + g_r$ - 1. Die Anzahl der Färbungen für G, die sich auf dem Rand von G unterscheiden, ergibt sich zu

$$2^{g_l-3} * 2^{g_r-3} * 3! > 2^{g_l+g_r-1+3} = 2^{g-3}$$

da wir drei Farben auf $G_l$ unabhängig von den Farben von $G_r$ permutieren können; eine Farbe ist durch die gemeinsame Ecke definiert.

Beispiel 2

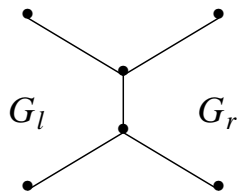

Es seien $G_l \in K$ und $G_r \in K$, mit zwei gemeinsamen Ecke bilden sie G (G ist kein i-Graph). Mit $g = \gamma(G)$, $g_l = \gamma(G_l)$, $g_r = \gamma(G_r)$ ist $g = g_l + g_r$ - 2. Die Anzahl der Färbungen ist analog

$$2^{g_l-3} * 2^{g_r-3} * 2 = 2^{g_l+g_r-2-3} = 2^{g-3}.$$

da wir zwei Farben auf $G_l$ unabhängig von den Farben von $G_r$ permutieren können. Zwei Farben sind durch die gemeinsamen Ecken definiert.

Die Graphen G in Beispiel 1 bzw. 2 haben eine trennende Eckenmenge [W, 72], mit 1 bzw. 2 Ecken, (deshalb kein i-Graph.) Die Färbungen unterscheiden sich auf dem Rand von G, sind aber nicht zu einem FS geordnet.

Weitere Beispiele

Beispiel 3

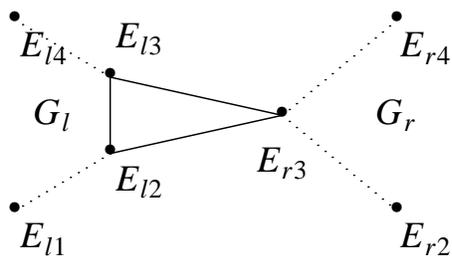

Es seien $G_l \in K$ und $G_r \in K$ mit FS, G wie abgebildet, G ist kein i-Graph. (Es ist nicht $G = G_l + G_r$ im Sinne unserer Definition für "+".) Mit $g = \gamma(G), g_l = \gamma(G_l), g_r = \gamma(G_r)$ ist $g = g_l + g_r$. Die Anzahl der Färbungen für G ergibt sich zu

$2^{g_l-3} * 2^{g_r-3} * 2 * 3! > 2^{g_l+g_3-6} * 2^3 = 2^{g-3}$ folgender maßen:

für $E_{r_3}$ sind zwei Farben möglich: $f(E_{r3}) \neq f(E_{l2})$, $f(E_{r3}) \neq f(E_{l3})$, die drei anderen Farben können wir auf $G_r$ permutieren unabhängig von den Färbungen von $G_l$.

Beispiel 4

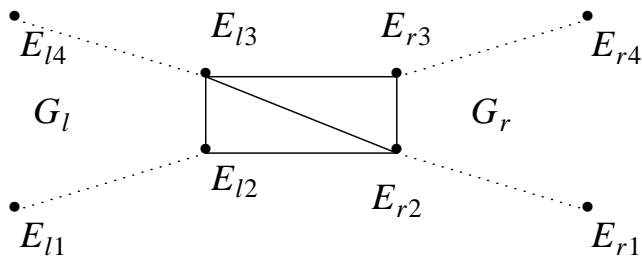

Es seien $G_l \in K$ und $G_r \in K$ mit FS, mit zwei gemeinsamen Ecke bilden sie G (G ist kein i-Graph). Mit $g = \gamma(G), g_l = \gamma(G_l), g_r = \gamma(G_r)$ ist $g = g_l + g_r$. Die Anzahl der Färbungen ist

$2^{g_l-3} * 2^{g_r-3} * 2 * 2 * 2 = 2^{g-3}$.

Für $E_{r_2}$ sind zwei Farben möglich ( $\neq f(E_{l2}), \neq f(E_{l3})$), dann sind für $E_{r3}$ wiederum zwei Farben möglich ( $\neq f(E_{r2}), \neq f(E_{l3})$), und dann können wir die zwei Farben verschieden von $f(E_{r2}), f(E_{r3})$ auf $E_r$ permutieren.

# Epilog

Diese Arbeit wurde bereits verschiedentlich veröffentlicht.

1. Im Jahre 1971 lag sie Prof. W. Jurkat vor während seines Aufenthaltes in Ulm, Abt. Mathematik I, Prof. A. Peyerimhoff.

2. Gedruckt als Technischer Bericht des Rechenzentrums der Universität Ulm.

3. Vortrag bei der Tagung Graphentheorie am Mathematischen Forschungsinstitut Oberwolfach im Juli 1975.

Ulm, August 2003

Lorenz Friess
Isnyerstr. 14
D-89079 Ulm
www.lorenz-friess.de